\documentstyle[amscd,12pt]{amsart}
\numberwithin{equation}{section}
\theoremstyle{plain}
   \newtheorem{thm}{Theorem}[section]
   \newtheorem{lem}[thm]{Lemma}
   \newtheorem{prop}[thm]{Proposition}
   \newtheorem{cor}[thm]{Corollary}

   \newtheorem{claim}[thm]{Claim}
\theoremstyle{definition}
   \newtheorem{defn}[thm]{Definition}
   
\theoremstyle{remark}
   \newtheorem{rem}[thm]{Remark}

\begin{document}

\title[$K_2$ of Tate curves]
{Surjectivity of $p$-adic regulator on $K_2$ of Tate curves}
\author[M. Asakura]
{Masanori Asakura}

\maketitle

%


\def\Spec{{\operatorname{Spec}}}
\def\Pic{{\mathrm{Pic}}}
\def\Ext{{\mathrm{Ext}}}
\def\NS{{\mathrm{NS}}}
\def\Picv{{\mathrm{Pic^0}}}
\def\Div{{\mathrm{Div}}}
\def\CH{{\mathrm{CH}}}
\def\deg{{\mathrm{deg}}}
\def\dim{{\operatorname{dim}}}
\def\codim{{\operatorname{codim}}}
\def\Coker{{\operatorname{Coker}}}
\def\ker{{\operatorname{ker}}}
\def\Image{{\operatorname{Image}}}
\def\Aut{{\mathrm{Aut}}}
\def\Hom{{\mathrm{Hom}}}
\def\Proj{{\mathrm{Proj}}}
\def\Sym{{\mathrm{Sym}}}
\def\Image{{\mathrm{Image}}}
\def\Gal{{\mathrm{Gal}}}
\def\GL{{\mathrm{GL}}}
\def\End{{\mathrm{End}}}
\def\P{{\bold P}}
\def\C{{\bold C}}
\def\R{{\bold R}}
\def\Q{{\bold Q}}
\def\Z{{\bold Z}}
\def\F{{\bold F}}
\def\l{\ell}
\def\ve{\varepsilon}
\def\lra{\longrightarrow}
\def\ra{\rightarrow}
\def\hra{\hookrightarrow}
\def\ot{\otimes}
\def\op{\oplus}
\def\vg{\varGamma}
\def\O{{\cal{O}}}
\def\ol#1{\overline{#1}}
\def\wt#1{\widetilde{#1}}
\def\us#1#2{\underset{#1}{#2}}
\def\os#1#2{\overset{#1}{#2}}
\def\lim#1{\us{#1}{\varinjlim}}
\def\plim#1{\us{#1}{\varprojlim}}

%


\def\mot{{\cal{M}}}
\def\Alb{{\mathrm{Alb}}}
\def\Perv{{\mathrm{Perv}}}
\def\gr{{\mathrm{gr}}}
\def\Gr{{\mathrm{Gr}}}
\def\rank{{\mathrm{rank}}}
\def\dlog{{\mathrm{dlog}}}
\def\Res{{\mathrm{Res}}}

\def\K{{\cal{K}}}
\def\D{{\Bbb{D}}}
\def\G{{\Bbb{G}}}
\def\dR{{\mathrm{dR}}}
\def\DR{{\mathrm{DR}}}
\def\Zar{{\mathrm{Zar}}}
\def\gen{{\mathrm{gen}}}
\def\rat{{\mathrm{rat}}}
\def\et{{\text{\'et}}}
\def\crys{{\mathrm{crys}}}
\def\Gr{{\mathrm{Gr}}}
\def\dego{\deg=0}
\def\Qb{\bar{\Q}}
\def\prim{0}
\def\al{\alpha}
\def\fp{{\frak{p}}}
\def\fm{{\frak{m}}}
\def\cExt{{\cal{E}}{\mathrm{xt}}}
\def\M{{\cal{M}}}
\def\D{\Delta}
\def\ord{{\mathrm{ord}}}
\def\reg{{\mathrm{reg}}}

\def\hsymbol{H_{K_2}}
\def\torsup{H^0_\Zar(E_K,{\cal K}_2)_{E[\infty]}}

\section{Introduction}
Let $X$ be a nonsingular projective variety over a field $K$ of
${\mathrm char}(K)\not=p$.
Then there are the {\it $p$-adic regulator maps}
\begin{equation}\label{padicintr}
c_{i,j}:K_i(X)\lra H^{2j-i}_\et(X,\Z_p(j)),
\quad i,~j\geq 0
\end{equation}
from Quillen's $K$-groups
to the \'etale cohomology groups with coefficients in the Tate twist
$\Z_p(j)$ (\cite{gillet}, \cite{soule2}).
When $K$ is a local field,
it is a long-standing problem whether
the maps $K_i(X)\ot\Q_p\to H^{2j-i}_\et(X,\Q_p(j))$ for $2j>i$
are surjective, in relation to the Beilinson conjectures (cf.
\cite{jan} \S 3).
The main result of this paper is to give an affirmative answer to
this problem for $K_2$ of
the Tate curves over certain $p$-adic fields:
\begin{thm}\label{mainthm}
Let $K$ be a finite extension of $\Q_p$.
Let $E_K=K^*/q^\Z$ be the Tate curve over $K$ where
$q\in K^*$ is a non-zero
element with its order $\ord (q)>0$.
Suppose that $K\subset \Q_p(\zeta)$ for some root of unity $\zeta$.
Then the $p$-adic regulator
\begin{equation}\label{qadicintr}
K_2(E_K)\ot\Q_p\lra H^2_\et(E_K,\Q_p(2))
\end{equation}
is surjective.
\end{thm}
The surjectivity is also true in the integral sense.
Namely due to Suslin's exact sequence (\cite{suslin} Cor.23.4),
Theorem \ref{mainthm} implies that
$H^0_\Zar(E_K,{\cal K}_2)\ot\Z_p
\to H^2_\et(E_K,\Z_p(2))$ is surjective, and it induces an isomorphism
$H^0_\Zar(E_K,{\cal K}_2)/p^\nu
\cong H^2_\et(E_K,\Z_p(2))/p^\nu$
for all $\nu\geq 1$ (Corollary \ref{pnuk2})

\medskip

We will prove Theorem \ref{mainthm} in the following steps.
Recall that there is a standard way to obtain elements of
$H^0_\Zar(E_K,{\cal K}_2))\ot\Q$ from torsion points of $E_K$
(e.g. \cite{deni-win} (5.1)).
The proof of Theorem \ref{mainthm}
is done by showing that such $K_2$-symbols
span the \'etale cohomology group $H^2_\et(E_K,\Q_p(2))$.
Due to the weight filtration on the \'etale cohomology,
$H^2_\et(E_K,\Q_p(2))$ is divided into
$H^1_\et(K,\Q_p(1))$ and $H^1_\et(K,\Q_p(2))$.
The proof of the first part (Part I, \S \ref{partipf})
is to show that $K_2$-symbols from torsion points
span $H^1_\et(K,\Q_p(1))$.
To do this, we will give a quite explicit formula
of the regulator maps. The technical results for it are given in
\S \ref{partipre}.
The second part (Part II, \S \ref{partiipf}) is to show that
some $K_2$-symbols from torsion points
span the image of the natural map
\begin{equation}\label{kctok}
\lim{F}H^1_{\et}(F,\Q_p(2)) \lra H^1_{\et}(K,\Q_p(2))
\end{equation}
where $F$ runs over all subfields of $K$
which are finite abelian extensions of $\Q$
(i.e. $F\subset K\cap
\Q(\mu)$ for some root $\mu$ of unity).
To do this, we will relate
some symbols in $K_2(E_K)$ with indecomposable elements of
$K_3(K)$ and apply a theorem of Soul\'e (Theorem \ref{k3et}).
\S \ref{partiipre} is the preliminary for it.
Finally, we will show that the map \eqref{kctok} is surjective
when $K\subset \Q_p(\zeta)$ for some root of unity $\zeta$
(Part III, \S \ref{galoissec}).
All over the proof (except Part III),
$p$-adic theta function is a basic tool.

\medskip

The surjectivity of \eqref{qadicintr} has an application to
torsion of $K_1(E_K)$ by Suslin's exact sequence:
\begin{thm}[Corollary \ref{torstr}]\label{main2}
Let the notations and assumption be as in Theorem \ref{mainthm}.
Then the torsion subgroup of $K_1(E_K)$ is finite.
More precisely,
let $\mu_n$ be the group of all roots of unity in $K$ where
$n$ denotes its cardinality.
Let $(~,~)_n:K^*/n\times K^*/n\to \mu_n$ be the Hilbert symbol
(cf. \cite{milnor} \S 15).
Then the torsion subgroup of $K_1(E_K)$ is isomorphic to
\begin{equation}\label{torsionpart}
\mu_n\op\mu_n\op\mu_n/(q,K^*)_n.
\end{equation}
\end{thm}
The decomposition in \eqref{torsionpart} corresponds to
the decomposition
$K_1(E_K)=K^*\op K^* \op V(E_K)$.

\medskip

There are previous works on the $l$-power torsion of $K_1$
for $l\not=p$.
T. Sato proved that
the $l$-adic regulator
$K_2(E_K)\ot\Q_l\ra H^2_\et(E_K,\Q_l(2))$ is surjective
and obtained the finiteness of the $l$-power torsion part of $K_1(E_K)$
(\cite{TSato}).
When $X$ is a nonsingular projective curve over a $p$-adic field
which has a good reduction,
the $l$-power torsion of $K_1(X)$ is finite and described by
the rational points of the jacobian of the special fiber
(Colliot-Th\'el\`ene and Raskind \cite{col-ras}).
See \S \ref{others} for details.
However, very little has been known about the $p$-power torsion
of $K_1$
or the surjectivity of $p$-adic regulators on $K_2$.

\medskip

Our proof of Theorem \ref{mainthm}
is comparable with T. Sato \cite{TSato}.
See \S \ref{tsato} (in particular Remark \ref{thesis}) for his proof.
However there is a big difference between $l$-adic and $p$-adic cases.
It is
based on the fact that $\dim H^2_\et(E_K,\Q_l(2))=1$ for $l\not=p$
whereas
$\dim H^2_\et(E_K,\Q_p(2))=2[K:\Q_p]+1$
where $[K:\Q_p]$ denotes the degree of $K$ over $\Q_p$.
In the $l$-adic case we only need to construct one $K_2$-symbol
which has nontrivial boundary (see \cite{scholl} for more calculation of
the boundary).
However it is not enough in the $p$-adic case. We need to calculate
$p$-adic regulators of symbols with trivial boundary.

\medskip

This paper is organized as follows.
\S \ref{preliminaries}
is the preliminaries on algebraic $K$-theory.
\S \ref{summary} is the summary of
Tate curves and theta functions.
\S \ref{partipre} -- \S \ref{galoissec} are devoted to prove
Theorem \ref{mainthm}. In \S \ref{applysec}, we give several
corollaries of Theorem \ref{mainthm}, including Theorem \ref{main2}.
In \S \ref{opentatesec}, we show that the $l$-adic regulator
on $K_2$ of any open subscheme of Tate curves is surjective.
\S \ref{opentatesec} is independent of the previous sections.

\medskip

\noindent{\it Acknowledgements}.
This paper was written during my stay
at the University of Chicago supported by
JSPS Postdoctoral Fellowships for Research Abroad.
I express sincere gratitude for their hospitality,
especially to professor Spencer Bloch.
I also thank professor Shuji Saito for sending me \cite{TSato},
and giving many valuable suggestions.

\section{Preliminaries on algebraic $K$-theory.}\label{preliminaries}
For an abelian group $M$,
we denote by $M[n]$ (resp. $M/n$)
the kernel (resp. cokernel) of multiplication by $n$.
\subsection{Higher $K$-theory and regulator maps}
Let $X$ be a separated quasi-projective scheme over a field $F$.
Let $P(X)$ be the exact category of locally free sheaves, and
$BQP(X)$ the simplicial set attached to $P(X)$ by Quillen (\cite{Q},
\cite{Sr}).
The {\it higher $K$-groups} of $X$ are defined as the homotopy
groups of $BQP(X)$:
$$
K_i(X)=\pi_{i+1}BQP(X), \quad i\geq 0.
$$
We refer \cite{Sr} for the general properties of higher $K$-theory
such as,
products, localization exact sequences, norm maps (also called transfer
maps) etc.

\medskip

Let $n$ be an integer which is prime to $\mathrm {char}(F)$.
Then there are the {\it regulator
maps}
\begin{equation}\label{chern}
c_{i,j}:K_i(X)\lra H^{2j-i}_\et(X,\Z/n(j)),
\quad i,~j\geq 0
\end{equation}
to the \'etale cohomology groups with coefficients in the Tate
twist $\Z/n(j)$ (\cite{gillet}, \cite{soule2}). They are
compatible with products, pull-backs and norm maps. When $X=\Spec
F$ and $i=j=1$, the regulator map is also known as the Galois
symbol
\begin{equation}\label{galoissymbol}
F^* \lra H^1_\et(F,\Z/n(1)), \quad f\longmapsto [f]
\end{equation}
in which $[f]$
is defined as the cocycle
$$
[f]:\Gal(\ol{F}/F) \lra \Z/n(1), \quad
\sigma \longmapsto \frac{\sigma(f^{1/n})}{f^{1/n}}.
$$

Of particular interest to us is the case that $i=j=2$ and $X$ is a curve.
\begin{lem}\label{comp11}
Let $X$ be a curve over $F$.
Put $X_{\ol{F}}=X\ot_F\ol{F}$.
Then the composition
\begin{equation}\label{comp428}
K_2(X)\lra H^2_\et(X,\Z/n(2))\lra
H^0(F,H^2_\et(X_{\ol{F}},\Z/n(2)))
\end{equation} is zero.
\end{lem}
\begin{pf}
We may assume $n=p^\nu$ with $p\not={\mathrm char}(F)$.
Since there is the isomorphism
$H^2_\et(X_{\ol{F}},\Z/p^\nu(2))\cong
\op H^2_\et(X_{i,\ol{F}},\Z/p^\nu(2))$ where $X_{i,\ol{F}}$ are
the irreducible components of $X_{\ol{F}}$,
we may assume that $X$ is irreducible.
If $X$ is not complete, there is nothing to prove because
of $H^2_\et(X_{\ol{F}},\Z/p^\nu(2))=0$.
Assume that $X$ is complete. Assume further
$H^0_\et(F,\Z_p(1))=0$.
Then, the assertion follows from the fact that
the composition \eqref{comp428} factors through
$$
\plim{\nu}
H^0(F,H^2_\et(X_{\ol{F}},\Z/p^\nu(2)))=H^0(F,\Z_p(1))=0.
$$
When $F$ is arbitrary, we choose an inductive limit $F=\varinjlim F_i$
and a projective limit $X=\varprojlim X_{F_i}$
where $F_i$ are finitely generated fields over the prime field
and $X_{F_i}$ are curves over $F_i$.
Since $H^0(F_i,\Z_p(1))=0$, we have
$$
K_2(X)=\lim{i}K_2(X_{F_i})\lra
\lim{i}H^2_\et(X_{\ol{F}_i},\Z/p^\nu(2)))=
H^2_\et(X_{\ol{F}},\Z/p^\nu(2)))
$$
is zero.
\end{pf}
By Lemma \ref{comp11} and the Hochschild-Serre spectral sequence,
the regulator map \eqref{chern} gives rise to a map
\begin{equation}\label{rho}
\rho_X:K_2(X)/n\lra
H^1(F,H^1_\et(X_{\ol{F}},\Z/n(2)))
\end{equation}
for a curve $X$.

\subsection{$K$-cohomology}\label{v(c)sec}
Let ${\cal K}_i$ be the Zariski
sheaf on $X$ associated to
the presheaf
$$
U\longmapsto K_i(U) \quad (U\subset X).
$$
Assume that $X$ is a nonsingular variety over $F$.
We denote by $X^i$ the set of points of height $i$.
Then the Gersten conjecture (proved by Quillen) says
the complex
$$
0\to {\cal K}_i \to K_i(F(X)) \to
\bigoplus_{x\in X^1} K_{i-1}(\kappa(x)) \to \cdots
\to \bigoplus_{x\in X^{\dim X}} K_{i-\dim X}(\kappa(x))
\to 0
$$
of Zariski sheaves is exact.
The above complex gives the flasque resolution of the sheaf
${\cal K}_i$. Therefore we have the isomorphism
\begin{equation}\label{coisom}
H^j_\Zar(X,{\cal K}_i)\cong
\frac{\ker (\bigoplus_{x\in X^j} K_{i-j}(\kappa(x)) \to
\bigoplus_{x\in X^{j+1}} K_{i-j-1}(\kappa(x)) )}{\Image
(\bigoplus_{x\in X^{j-1}} K_{i-j+1}(\kappa(x)) \to
\bigoplus_{x\in X^{j}} K_{i-j}(\kappa(x)) )}
\end{equation}
In particular,
when $X$ is a nonsingular curve over $F$, we have
the exact sequence
\begin{equation}\label{coisomcurve}
0\ra H^0_\Zar(X,{\cal K}_2)
\ra K_2^M(F(X)) \os{\tau}{\ra}
\bigoplus_{x\in X^1} \kappa(x)^*
\ra H^1_\Zar(X,{\cal K}_2) \ra 0.
\end{equation}
Here $K_2^M$ denotes the Milnor $K$-theory (which coincides with
Quillen's $K_2$ by Matsumoto's theorem), and $\tau=\sum \tau_x$ is
the sum of the tame symbol $\tau_x$ at $x\in X^1$:
\begin{equation}\label{tame}
\tau_x:K_2^M(F(X))\lra \kappa(x),\quad
\{f,g\}\mapsto(-1)^{{\mathrm ord}_x(f)
{\mathrm ord}_x(g)}\frac{f^{{\mathrm ord}_x(g)}}
{g^{{\mathrm ord}_x(f)}}.
\end{equation}
Hereafter,
we always identify the $K$-cohomology $H^0_\Zar(X,{\cal K}_2)$
(resp. $H^1_\Zar(X,{\cal K}_2)$) with the kernel of $\tau$
(resp. cokernel of $\tau$) for a nonsingular curve $X$.
Due to the localization exact sequence of $K$-theory, we see that
there is a natural surjection
$K_2(X)\to H^0_\Zar(X,{\cal K}_2)$ and
an exact sequence
\begin{equation}\label{k1coh}
0\lra H^1_\Zar(X,{\cal K}_2) \lra K_1(X) \lra F(X)^*
\os{\mathrm ord}{\lra}
\bigoplus_{x\in X^1} \Z.
\end{equation}
Suppose further that $X$ is a complete nonsingular curve.
Then the norm maps $N_{\kappa(x)/F}:\kappa(x)^*\to F^*$
induce the norm map
$H^1_\Zar(X,{\cal K}_2)\to F^*$ on $K$-cohomology.
We denote by $V(X)$ the kernel of it:
\begin{equation}\label{v(c)}
0\lra V(X)\lra H^1_\Zar(X,{\cal K}_2) \lra F^*.
\end{equation}
If $X$ has a $F$-rational point
the right map is surjective and we have a decomposition
$$
K_1(X)=F^*\op H^1_\Zar(X,{\cal K}_2)
=F^*\op F^*\op V(X).
$$

\subsection{Suslin's exact sequence.}\label{suslinseq}
Let $X$ be a nonsingular curve over $F$.
It follows from the Riemann-Roch theorem (\cite{gillet}) that
the regulator map $c_{2,2}:K_2(X)\to H^2_\et(X,\Z/n(2))$
induces a map
$H^0_\Zar(X,{\cal K}_2)
\to H^2_\et(X,\Z/n(2))
$.
A.Suslin proved that
there is the natural exact sequence
(\cite{suslin} Cor.23.4)
\begin{equation}\label{universalc}
0\lra H^0_\Zar(X,{\cal K}_2)/n
\lra H^2_\et(X,\Z/n(2)) \lra
H^1_\Zar(X,{\cal K}_2)[n]\lra 0
\end{equation}
for ${\mathrm char}(F)\not\vert n$. (It is proved not only for
curves but also for any nonsingular varieties. However, it is not
used in this paper.) Suslin's sequence \eqref{universalc} will be
used for the proof of Theorem \ref{main2}.

By Lemma \ref{comp11} and the Hochschild-Serre spectral sequence,
we have a map
\begin{equation}\label{rhoz}
H^0_\Zar(X,{\cal K}_2)/n\lra
H^1(F,H^1_\et(X_{\ol{F}},\Z/n(2))),
\end{equation}
which is compatible with \eqref{rho} under the natural
surjection $K_2(X) \to H^0_\Zar(X,{\cal K}_2)$.
Without confusing, we also write the map \eqref{rhoz} by $\rho_X$.


\section{Tate curves and $p$-adic theta functions}\label{summary}
\def\RR{R}
\def\KK{K}

We give a brief review on Tate curves and theta functions.
No proofs are in this section.
A good reference is Silverman's book \cite{silverman}.
\subsection{}
Let $K$ be a finite extension of $\Q_p$, and
${\mathrm ord}_K:K^*\to \Z$ the map of order
such that ${\mathrm ord}_K(\pi_K)=1$
where $\pi_K$ denotes a uniformizer of $K$.
Let $q\in K^*$ satisfy ${\mathrm ord}_K(q)>0$.
The {\it Tate curve} $E_K$ with the period $q$ is defined
as the elliptic curve over $K$ defined by the equation
\begin{equation}\label{tate0}
y^2+xy=x^3+a_4(q)x+a_6(q)
\end{equation}
where
\begin{equation}\label{tate1}
a_4(q)=-5\sum_{n=1}^\infty \frac{n^3q^n}{1-q^n},\quad
a_6(q)=-\sum_{n=1}^\infty \frac{(5n^3+7n^5)q^n}{1-q^n}.
\end{equation}
This is a $p$-adic analogue of the complex torus $\C^*/q^\Z$. As
is so in the classical case, the discriminant $\Delta$ of $E_K$ is
given by
$$
\Delta=q\prod_{n\geq 1}(1-q^n)^{24}
$$
and the $j$-invariant
$$
j(E_K)=\frac{1}{q}+744+196884q+\cdots.
$$
The series
$$
X(u)=\sum_{n\in\Z}\frac{q^nu}{(1-q^nu)^2}
-2\sum_{n\geq1}\frac{nq^n}{1-q^n}
$$
$$
Y(u)=\sum_{n\in\Z}\frac{(q^nu)^2}{(1-q^nu)^3}
+\sum_{n\geq1}\frac{nq^n}{1-q^n}
$$
converge for all $u\in \ol{K}-q^\Z$. They induce a bijective
homomorphism
\begin{equation}\label{thehom}
\ol{K}^*/q^\Z\os{\cong}{\lra} E_K(\ol{K}),
\quad
u\longmapsto
\begin{cases}
(X(u),Y(u)) & \text{if }u\not\in q^\Z\\
O & \text{if }u\in q^\Z
\end{cases}
\end{equation}
where $O\in E_K(K)$ denotes the infinity point.
We often identify $E_K(\ol{K})$ with $\ol{K}^*/q^\Z$ by the
isomorphism \eqref{thehom}.

\begin{defn}[Theta function]\label{thetadefn}
$$
\theta(u)=\theta(u,q)
\os{{\text{def}}}{=}(1-u)\prod_{n=1}^\infty(1-q^nu)(1-q^nu^{-1})
$$
\end{defn}
$\theta(u)$ converges for all $u\in \ol{K}^*$ and satisfy
\begin{equation}\label{thetaformula1}
\theta(qu)=\theta(u^{-1})=-\frac{1}{u}\theta(u).
\end{equation}
Using \eqref{thetaformula1}, we see that
a function
\begin{equation}\label{11}
f(u)=c\prod_i\frac{\theta(\alpha_iu)}{\theta(\beta_iu)}
\end{equation}
is $q$-periodic if $\prod_i{\alpha_i}/{\beta_i}=1$.
Conversely, for any rational function $f(u)$ on $E_{\ol{K}}:=E_K\ot
\ol{K}$,
one can find $c,~\alpha_i,~\beta_i\in\ol{\KK}^*$ such that $f(u)$
is given as in \eqref{11}. Thus we have a one-one correspondence
\begin{equation}\label{cor}
\ol{\KK}(E_{\ol{\KK}})^*
\os{1:1}{\longleftrightarrow}
\left\{
c\prod_i\frac{\theta(\alpha_iu)}{\theta(\beta_iu)}
~\vert~
c,~\alpha_i,~\beta_i\in\ol{\KK}^*
\text{ with }\prod_i{\alpha_i}/{\beta_i}=1
\right\}.
\end{equation}
We often identify the both sides of \eqref{cor}.
Since the correspondence \eqref{cor} is compatible with the action
of the Galois group $G_K$, a rational function
$f(u)\in \ol{K}(E_{\ol{K}})^*$ is contained
in $K(E_K)$ if and only if $c\in K^*$ and
\begin{equation}\label{divisor1}
\sum_i[\alpha_i]-[\beta_i]=\sum_i[\alpha^\sigma_i]
-[\beta^\sigma_i]
\end{equation}
as divisor on ${\Bbb G}_{m,\ol{K}}$ for all
$\sigma\in {\mathrm Gal}(\ol{K}/K)$.

\subsection{Definition of $K\langle u\rangle$}
Let $R$ be the integer ring of $K$.
We define a ring
\begin{align*}
\RR\langle u\rangle&\os{\mathrm def}{=}\plim{n}\RR/\pi_K^n[[u]][u^{-1}]\\
&=\{\sum_{i=-\infty}^{+\infty}a_iu^i\in R[[u,u^{-1}]]~;~
\ord_K(a_i)\to+\infty \text{ as }i\to-\infty\}.
\end{align*}
This is a discrete valuation ring with
a uniformizer $\pi_K$.
We write by $K\langle u\rangle$ the quotient field of $\RR\langle u\rangle$:
$$\KK\langle u\rangle=\RR\langle u\rangle[\pi_K^{-1}]=
\RR\langle u\rangle\ot_R\KK.$$
The field $K\langle u\rangle$ contains $K(u)$, but not $K((u))$.

Since the theta function $\theta(u)$ is contained in
$R\langle u\rangle$, the correspondence \eqref{cor}
defines an inclusion $K(E_K) \hookrightarrow
K\langle u\rangle$. Thus we have a dominant morphism
\begin{equation}\label{cor2}
\Spec\KK\langle u\rangle\lra E_\KK.
\end{equation}
of schemes.
\subsection{Semistable reduction of Tate curves}\label{reductionsec}
Let ${\cal C}/R$ be a minimal proper regular model of $E_K/K$
over the integer ring $R$.
Put $n={\mathrm ord}_K(q)$.
By a result of Kodaira and N\'eron, the special fiber
$Y=Y_1+\cdots+Y_n$ is type $I_n$, namely, if $n\geq 2$,
$Y_i$ are nonsingular rational curves which are arranged
in the shape of a $n$-gon, and if $n=1$, $Y=Y_1$ is an irreducible
rational curve with a node.

Let ${\cal E}/R$ be the N\'eron model of $E_K/K$. It is the
largest subscheme of ${\cal C}/R$ which is smooth
(\cite{silverman} Theorem 6.1). The group law on $E_K$ extends to
make ${\cal E}/R$ into a commutative group scheme over $R$. The
special fiber $E_0$ of ${\cal E}/R$ is a commutative group scheme
over the residue field $k$ which consists of $n$-copies of ${\Bbb
G}_{m,k}$. More precisely, we have an isomorphism $E_0\cong {\Bbb
G}_{m,k}\times \Z/n$ as group schemes. By the N\'eron mapping
property, we have ${\cal E}(R)=E_K(K)$. Therefore we have a
homomorphism $E_K(K)=K^*/q^\Z\to E_0(k)=k^*\times \Z/n$. It is
explicitly given by $aq^{i/n}\mapsto (a~{\mathrm
mod}\pi_K,i~{\mathrm mod}n)$.

For an integer $r\geq 1$,
we put $R_r:=R/\pi_K^{r+1}$ and ${\cal E}_r:= {\cal E}\ot_RR_r$.
Then ${\cal E}_r/R_r$ is a group scheme.
Let ${\cal E}_r^o$ be the identity component of ${\cal E}_r$.
As we have seen in the above, ${\cal E}_0^o\cong {\Bbb G}_{m,k}$.
Due to the rigidity of algebraic tori (\cite{SGA3} exp.IX \S3),
there is an isomorphism
${\cal E}^o_r\cong{\Bbb G}_{m,R_r}$ of group schemes over $R_r$.
The embedding
$$
h:{\Bbb G}_{m,R_r} \lra {\cal E}_r\subset {\cal E}\subset {\cal C}
$$
is (locally) defined by
$$
h^*(x)=X(u)~{\mathrm mod}\pi_K^{r+1},\quad
h^*(y)=Y(u)~{\mathrm mod}\pi_K^{r+1}.
$$
Taking the inductive limit, we have a homomorphism
$$
{\Bbb G}_{m}^{\mathrm for}:=
\lim{r}{\Bbb G}_{m,R_r}\lra {\cal E}^{\mathrm for}:=
\lim{r}{\cal E}_r
$$
of formal schemes. (Note that it is not algebraizable.)
Composing with ${\cal E}^{\mathrm for}\to {\cal C}$,
we get a morphism
${\Bbb G}_{m}^{\mathrm for}\ot_R K\to {\cal C}\ot_RK=E_K$.
It gives a homomorphism
$$
K(E_K) \lra \left(\plim{r}R/\pi_K^r[u,u^{-1}]\right)\ot_RK\lra
\left(\plim{r}R/\pi_K^r[[u]][u^{-1}]\right)\ot_RK
=K\langle u\rangle
$$
of fields. This gives another definition of \eqref{cor2}.


\section{The weight exact sequence}\label{partipre}
\subsection{Weight exact sequence.}\label{weightseq}
The algebraic fundamental group $\pi_1(E_{\ol{\KK}})$
of the Tate curve
$E_{\ol{\KK}}:=E_K\ot\ol{K}$ is isomorphic to $\hat{\Z}\times\hat{\Z}$
since the characteristic of $\KK$ is zero.
Let us give its generators explicitly.
For a Galois covering $f:X\to  E_{\ol{\KK}}$, we denote by
$\Aut(f)$ the group of $\ol{\KK}$-automorphisms
$T:E_{\ol{\KK}}\to E_{\ol{\KK}}$ such that $fT=f$.
Let $\nu_n:E_{\ol{\KK}}\to E_{\ol{\KK}}$ be the Galois
covering given by the multiplication $x\mapsto x^n$.
Then $\Aut(\nu_n)$ is isomorphic to $\Z/n\times \Z/n$.
The generators are translations
\begin{equation}\label{trans1}
T_{\zeta_n}:\ol{\KK}^*/q^\Z\lra \ol{\KK}^*/q^\Z, \quad x\to x\zeta_n ,
\end{equation}
\begin{equation}\label{trans2}
T_{q^{1/n}}:\ol{\KK}^*/q^\Z\lra \ol{\KK}^*/q^\Z, \quad x\to xq^{1/n}
\end{equation}
where $\zeta_n$ is a primitive $n$-th root of unity.
$\pi_1(E_{\ol{\KK}})$ is isomorphic to $\plim{n}\Aut(\nu_n)$,
and its (topological) generators are $T_{\zeta_\infty}=\plim{n}T_{\zeta_n}$
and $T_{q_\infty}=\plim{n}T_{q^{1/n}}$.
There is the fibration exact sequence
\begin{equation}\label{fib0}
1\lra \pi_1(E_{\ol{\KK}})
\lra
\pi_1(E_\KK)
\lra G_\KK \lra 1
\end{equation}
where $G_\KK$ denotes the absolute Galois group of $\KK$.
The sequence \eqref{fib0} is split by the map coming from a
$\KK$-rational point
$\Spec\KK\to E_\KK$. Thus $\pi_1(E_\KK)$ is isomorphic to the semidirect
product $\pi_1(E_{\ol{\KK}})\cdot G_\KK$ with
\begin{equation}\label{semidirect}
\sigma T_{\zeta_\infty}\sigma^{-1}=T_{\zeta^\sigma_\infty}
=\plim{n}T_{\zeta^\sigma_n},
\quad
\sigma T_{q_\infty}\sigma^{-1}
=T_{q^\sigma_\infty}=\plim{n}T_{(q^{1/n})^\sigma}
\end{equation}
for $\sigma\in G_\KK$.
Denote by $\langle T_{\zeta_\infty}\rangle
\subset \pi_1(E_{\ol{\KK}})$ the closed subgroup
generated by $T_{\zeta_\infty}$, and $\pi\subset\pi_1(E_\KK)$
the closed subgroup generated by $\langle T_{\zeta_\infty}\rangle$
and $G_\KK$:
$$
\begin{CD}
1@>>> \pi_1(E_{\ol{\KK}})
@>>>
\pi_1(E_\KK)
@>>> G_\KK
@>>> 1\\
@.@AAA@AAA@AA{=}A\\
1@>>> \langle T_{\zeta_\infty}\rangle
@>>>
\pi
@>>> G_\KK
@>>> 1.
\end{CD}
$$
Due to \eqref{semidirect}, $\pi$ is isomorphic to the semidirect product
$\langle T_{\zeta_\infty}\rangle\cdot G_\KK$ with the relation
$\sigma T_{\zeta_\infty}\sigma^{-1}=T_{\zeta^\sigma_\infty}$.
Therefore the natural map
$$
H^1_\et(E_{\ol{\KK}},\Z/n(j+1))=
\Hom(\pi_1(E_{\ol{\KK}}),\Z/n(j+1))\lra
\Hom(\langle T_{\zeta_\infty}\rangle,\Z/n(j+1))
$$
is compatible with respect to
$G_\KK$-action, and the target is
isomorphic to $\Z/n(j)$ as $G_\KK$-module.
Similarly we can see that the map
$$
\Hom(\pi_1(E_{\ol{\KK}})/\langle T_{\zeta_\infty}\rangle
,\Z/n(j+1))\lra
\Hom(\pi_1(E_{\ol{\KK}}),\Z/n(j+1))
$$
is compatible with $G_\KK$-action, and
the source is
isomorphic to $\Z/n(j+1)$ as $G_\KK$-module.
As a result, we have an exact sequence of $G_\KK$-modules:
\begin{equation}\label{van0}
0\lra \Z/n(j+1)\lra
H^1_\et(E_{\ol{\KK}},\Z/n(j+1))
\lra
\Z/n(j) \lra 0.
\end{equation}
This is called the {\it weight exact sequence}.

\begin{lem}\label{psilem}
Let $E_{m,K}=K^*/q^{m\Z}$.
Denote by $\psi_m:E_K\to E_{m,K}$ and $\phi_m:E_{m,K}\to E_K$
the homomorphism given by $x\mapsto x^m$ and
the natural surjection respectively.
Then the pull-backs $\psi^*_m$ and $\phi^*_m$
induce the following commutative
diagrams:
$$
\begin{CD}
0@>>> \Z/n(j) @>>> H^1_\et(E_{m,\ol{K}},\Z/n(j))
@>>> \Z/n(j-1) @>>>0\\
@.@V{=}VV@V{\psi_m^*}VV@VV{{\text {\rm mult. by } } m}V\\
0@>>> \Z/n(j) @>>> H^1_\et(E_{\ol{K}},\Z/n(j))
@>>> \Z/n(j-1) @>>>0,
\end{CD}
$$
$$
\begin{CD}
0@>>> \Z/n(j) @>>> H^1_\et(E_{\ol{K}},\Z/n(j))
@>>> \Z/n(j-1) @>>>0\\
@.@V{{\text {\rm mult. by } } m}VV@V{\phi_m^*}VV@VV{=}V\\
0@>>> \Z/n(j) @>>> H^1_\et(E_{m,\ol{K}},\Z/n(j))
@>>> \Z/n(j-1) @>>>0.
\end{CD}
$$
\end{lem}
\begin{pf}
The map
$\psi_{m*}:\pi_1(E_K)\to \pi_1(E_{m,K})$ is given as follows
$$
\psi_{m*}T_{\zeta_\infty}=\plim{i}\psi_{m*}T_{\zeta_i}
=\plim{i}T_{\zeta^m_i}=T_{\zeta_\infty}^m,
$$
$$
\psi_{m*}T_{q_\infty}=\plim{i}\psi_{m*}T_{q^{1/i}}=
\plim{i}T_{q^{m/i}}\equiv \plim{i}T_{(q^m)^{1/i}}
=T_{(q^m)_\infty} \mod \langle T_{\zeta_\infty}\rangle.
$$
Thus the commutative diagram for $\psi_m$ follows.
The diagram for $\phi_m$ follows in a similar way.
\end{pf}

\subsection{Definition of $\tau^\et_\infty$}\label{rwsec}
From the weight exact sequence \eqref{van0},
we have
$$
\begin{CD}
@. H^0_\Zar(E_\KK,{\cal K}_2)/n@.\\
@.@V{\rho}VV@.\\
H^1(\KK,\Z/n(2))
@>{a}>>
H^1(\KK,H^1_\et(E_{\ol{\KK}},\Z/n(2)))
@>{b}>>
H^1(\KK,\Z/n(1))=\KK^*/n.
\end{CD}
$$
Here $\rho=\rho_{E_K}$ is as in \eqref{rhoz}.
We define $\tau_\infty$ as the composition of $b$ and $\rho$:
$$
\tau^\et_\infty\os{\text{def}}{=}b\cdot\rho:
H^0_\Zar(E_\KK,{\cal K}_2)/n\lra
\KK^*/n.
$$
By the construction,
the maps $\rho$ and $\tau^\et_\infty$ are
compatible with the pull-back and the norm map
for any finite extension $L/\KK$.

\subsection{}

\def\KKdag{{\KK\langle u\rangle_{\ol{\KK}}}}

Put $\KKdag= \KK\langle u\rangle\ot_\KK\ol{\KK}$. Consider a map
\begin{equation}\label{res}
(\KKdag)^*/n
\lra \Z/n, \quad
f\longmapsto {\mathrm Res}\frac{df}{f}
\end{equation}
where ${\mathrm Res}$ denotes the residue map at $u=0$, namely if
we express $\omega=\sum_{n\in\Z}a_nu^ndu$ in the unique way then
${\mathrm Res}(\omega)=a_{-1}$. The map \eqref{res} is clearly a
homomorphism of $G_\KK$-module. On the other hand the morphism
\eqref{cor2} induces
\begin{equation}\label{h1}
H^1_\et(E_{\ol{\KK}},\Z/n(1))
\lra
H^1_\et(\KKdag,\Z/n(1))
=
(\KKdag)^*/n.
\end{equation}

\begin{lem}\label{thedilem}
The diagram
\begin{equation}\label{thedi}
\begin{CD}
H^1_\et(E_{\ol{\KK}},\Z/n(1))@>>> \Z/n\\
@VVV@VV{=}V\\
(\KKdag)^*/n
@>>> \Z/n
\end{CD}
\end{equation}
is commutative. Here the maps are as in \eqref{van0},
\eqref{res} and \eqref{h1}.
\end{lem}
\begin{pf}
Fix $q^{1/n}\in\ol{\KK}$ and a primitive $n$-th root of unity
$\zeta_{n}$.
We put
$$
f_1(u):=\frac{\theta(q^{1/n}u)^{n}}{\theta(u)^{n-1}\theta(qu)}
=-u\left(\frac{\theta(q^{1/n}u)}{\theta(u)}
\right)^{n}
$$
$$
f_2(u):=\left(\frac{\theta(\zeta_{n}u)}{\theta(u)}
\right)^{n}.
$$
The divisors of $f_1$ and $f_2$ are $n([q^{-1/n}]-[1])$
and $n([\zeta_{n}^{-1}]-[1])$ respectively where
$[\alpha]$ denotes the divisor of a closed point $\alpha\in\ol{K}^*/q^\Z$.
Therefore, each $f_i$ defines the cohomology class
$[f_i]\in H^1_\et(E_{\ol{\KK}},\Z/n(1))$.
We claim that $[f_1]$ and $[f_2]$ span the cohomology group
$H^1_\et(E_{\ol{\KK}},\Z/n(1))$.
Recall that the cohomology class $[f_i]$ is defined as
$$
T_{\zeta_\infty}\mapsto T_{\zeta_{n}}(\phi_{n}^*f_i^{1/n})/
\phi_{n}^*f_i^{1/n},
\quad
T_{q_\infty}\mapsto T_{q^{1/n}}(\phi_{n}^*f_i^{1/n})/
\phi_{n}^*f_i^{1/n}
$$
under the isomorphism $H^1_\et(E_{\ol{\KK}},\Z/n(1))
\cong\Hom(\pi_1(E_{\ol{\KK}}),\Z/n(1))$.
Note
$$
\phi_{n}^*f_1^{1/n}
=(-1)^{1/n}v\frac{\theta(q^{1/n}v^{n})}{\theta(v^{n})},
\quad
\phi_{n}^*f_2^{1/n}=
\frac{\theta(\zeta_{n}v^{n})}{\theta(v^{n})}
$$
and $T_{\zeta_{n}}$ and $T_{q^{1/n}}$ are given by
$v\mapsto \zeta_{n}v$ and $v\mapsto q^{1/n}v$ respectively.
Therefore we see
$$
[f_1]:
\begin{matrix}
T_{\zeta_\infty} & \mapsto &\zeta_{n}\\
T_{q_\infty} & \mapsto &0
\end{matrix},
\quad
[f_2]:
\begin{matrix}
T_{\zeta_\infty} & \mapsto &0\\
T_{q_\infty} & \mapsto & \zeta_{n}^{-1}.
\end{matrix}
$$
This shows that $[f_1]$ and $[f_2]$ span
$H^1_\et(E_{\ol{\KK}},\Z/n(1))=
\Hom(\pi_1(E_{\ol{\KK}}),\Z/n(1))$.

To show the commutativity of the diagram \eqref{thedi},
it suffices to show that
$$
{\mathrm Res}\frac{df_1}{f_1}=1,\quad
{\mathrm Res}\frac{df_2}{f_2}=0
\quad\mod n\Z.
$$
Each of them is straightforward.
\end{pf}

Put $\left(K_2(\KK\langle u\rangle)/n\right)'$ by
$$
0\lra \left(K_2(\KK\langle u\rangle)/n\right)'
\lra K_2(\KK\langle u\rangle)/n
\lra H^2_\et(\KKdag,\Z/n(2)).
$$
Then we have a commutative diagram
$$
\begin{CD}
K_2(E_\KK)/n @>{\rho}>> H^1(K,H^1_\et(E_{\ol{\KK}},\Z/n(2)))
@>{b}>>K^*/n\\
@VVV@VVV@VV{=}V\\
\left(K_2(\KK\langle u\rangle)/n\right)'
@>>> H^1_\et(K,H^1_\et(\KKdag,\Z/n(2)))
@>>>K^*/n
\end{CD}
$$
where the commutativity of the right square is due to
Lemma \ref{thedilem}. We denote
by $\hat{\tau}^\et_\infty$
the composition of the below arrows:
\begin{equation}\label{hochserre2}
\hat{\tau}^\et_\infty:
\left(K_2(\KK\langle u\rangle)/n\right)'
\lra
\KK^*/n.
\end{equation}

\begin{lem}
Let $\KKdag^*/n\to\Z/n$ be as in \eqref{res}.
Then the
composition $\ol{\KK}(u)^*/n\to\KKdag^*/n\to\Z/n$
is given by
$$
f\longmapsto \sum_{\alpha}{\mathrm Res}_\alpha \frac{df}{f}
$$
where $\alpha$ runs over all $\alpha\in \ol{\KK}$ such that
${\mathrm ord}_K(\alpha)>0$ (including $\alpha=0$).
\end{lem}
\begin{pf}
This is straightforward because
$$
(u-\alpha)^{-1}=
\begin{cases}
u^{-1}\sum_{n=0}^\infty(\alpha u^{-1})^n & {\mathrm ord}_K(\alpha)>0\\
-\alpha^{-1}\sum_{n=0}^\infty(\alpha^{-1} u)^n
& {\mathrm ord}_K(\alpha)\leq 0
\end{cases}
$$
in $\KKdag$.
\end{pf}
It is well-known that
the composition
$$
K_2(\KK(u))/n \lra H^1(\KK,H^1_\et(\ol{\KK}(u),\Z/n(2)))
\os{{\mathrm Res}_\alpha}{\lra}H^1(\KK,\Z/n(1))=\KK^*/n
$$
coincides with the tame symbol $\tau_\alpha$ at $u=\alpha$.
Therefore the following map
\begin{equation}\label{a1}
K_2(\KK(u))/n \lra
\left(K_2(\KK\langle u\rangle)/n\right)'
\os{\hat{\tau}^\et_\infty}{\lra}
\KK^*/n.
\end{equation}
is given by
$$
\xi \longmapsto
\sum_{\alpha}\tau_{\alpha} (\xi)
$$
where $\alpha$ runs over all $\alpha\in \ol{\KK}$ such that
${\mathrm ord}_K(\alpha)>0$ (including $\alpha=0$).

\medskip

Summarizing the above results, we have the following:
\begin{thm}\label{formula}
The diagram
$$
\begin{CD}
K_2(E_\KK)/n
@>{\tau^\et_\infty}>>
\KK^*/n\\
@VVV@VV{=}V\\
\left(K_2(\KK\langle u\rangle)/n\right)'
@>{\hat{\tau}^\et_\infty}>>
\KK^*/n\\
@AAA@AA{=}A\\
K_2(\KK(u))/n
@>{\sum_{\alpha}\tau_{\alpha}}>>
\KK^*/n
\end{CD}
$$
is commutative
where $\alpha$ runs over all $\alpha\in \ol{\KK}$ such that
${\mathrm ord}_K(\alpha)>0$ (including $\alpha=0$).
\end{thm}
This theorem enables us to calculate the map $\tau^\et_\infty$ explicitly
(cf. proof of Proposition \ref{sa}).


\section{Proof of Theorem \ref{mainthm} : Part I}\label{partipf}
There is the Hochschild-Serre spectral sequence
$$
E_2^{ij}=H^i_\et(K,H^j(E_{\ol{K}},\Q_p(2)))\Longrightarrow
H_\et^{i+j}(E_K,\Q_p(2)).
$$
It degenerates at $E_2$-terms.
Since the cohomological dimension of $K$ is 2 (\cite{galois} II 4.3),
we have $E_2^{ij}=0$ for $i\geq 3$.
Moreover, since $[K:\Q_p]<\infty$,
$$
E_2^{02}=H^0_\et(K,H^2(E_{\ol{K}},\Q_p(2)))=H^0_\et(K,\Q_p(1))=0.
$$
Due to the duality theorem for the Galois cohomology
of local fields (loc.cit. II. 5.2. Theorem. 2), we have
$$
H^2_\et(K,\Z/p^\nu(2))\cong
\Hom(H^0_\et(K,\Z/p^\nu(-1)),\Q/\Z)\cong H^0_\et(K,\Z/p^\nu(1)),
$$
and hence $E_2^{20}=H^2_\et(K,\Q_p(2))=0$.
Therefore we have an isomorphism
\begin{equation}\label{open1}
H_\et^{2}(E_K,\Q_p(2))\cong H^1_\et(K,H^1(E_{\ol{K}},\Q_p(2))).
\end{equation}
Thus in order to prove Theorem \ref{mainthm}
it suffices to prove that the cardinality of the cokernel of the map
\begin{equation}\label{qpmap0}
\rho:H^0_\Zar(E_K,{\cal K}_2)/p^\nu
\lra H^1_\et(K,H^1(E_{\ol{K}},\Z/p^\nu(2)))
\end{equation}
has an upper bound which
does not depend on $\nu$.
Due to the weight exact sequence \eqref{van0}, we have
an exact sequence
$$
H^1_\et(K,\Z/p^\nu(2))
\os{a}{\lra}
H^1_\et(K,H^1_\et(E_{\ol{K}},\Z/p^\nu(2)))
\os{b}{\lra}
H^1_\et(K,\Z/p^\nu(1))=K^*/p^\nu.
$$
Note that both of the kernel of $a$ and the cokernel of $b$
are finite whose orders are at most
$\sharp K^*[p^\infty]$.
Then, we first show the following.
\begin{description}
\item[(Part I)]
The cardinality of the cokernel of the map
$$
\tau^\et_\infty:H^0_\Zar(E_K,{\cal K}_2)/p^\nu\lra K^*/p^\nu
$$
has an upper bound which does not depend
on $\nu$.
\end{description}

\medskip
Second we put
\begin{equation}\label{mot}
\hsymbol:=\Image (\ker ~\tau^\et_\infty\to H^1_\et(K,\Z/p^\nu(2))/\ker~a)
\end{equation}
and
\begin{equation}\label{ab}
H_{\mathrm ab}:=\Image
(\lim{F}
H^1_\et(F,\Z_p(2))
\to H^1_\et(K,\Z/p^\nu(2))/\ker ~a),
\end{equation}
where $F$ runs over all subfields of $K$
which are finite abelian extensions of $\Q$.
\begin{description}
\item[(Part II)]
$\hsymbol
\supset mH_{\mathrm ab}$
for some $m\not=0$ which does not depend on $\nu$.
\end{description}
Final step is to show that the index of $H_{\mathrm ab}$
has an upper bound which
does not depend on $\nu$, or
equivalently
\begin{description}
\item[(Part III)]
The map \eqref{kctok}
is surjective if $K\subset \Q(\zeta)$ for some root of unity $\zeta$.
\end{description}

\begin{rem}
We do not need any assumption on $K$ for the proofs of
(Part I) and (Part II).
Therefore, we have the surjectivity of
the $p$-adic regulator \eqref{qadicintr}
only if $K$ satisfies that \eqref{kctok} is surjective.
\end{rem}

\subsection{Proof of Part I : Step 1}\label{tsato}
Let ${\mathrm ord}_K:K^*\to \Z$ be the map of order such that
${\mathrm ord}_K(\pi_K)=1$ for a uniformizer $\pi_K\in K$.
We first show that
the map
\begin{equation}\label{sato}
{\mathrm ord}_K\cdot\tau^\et_\infty:H^0_\Zar(E_K,{\cal K}_2)/p^\nu
\lra
\Z/p^\nu
\end{equation}
is surjective.
More precisely, let
\begin{equation}\label{sato2}
o_K:H^0_\Zar(E_K,{\cal K}_2) \lra
\plim{n}H^0_\Zar(E_K,{\cal K}_2)/n
\lra
\plim{n}\Z/n=
\hat{\Z}
\end{equation}
be the composition.
Then we construct a symbol $\xi\in H^0_\Zar(E_K,{\cal K}_2)$
(which comes from torsion points of $E_K$)
such that
$o_K(\xi)$ is a nonzero integer.

\medskip

Let $L/K$ be a finite extension such that there is a uniformizer $\pi_0$
of $L$ satisfying $q=\pi_0^r$ for some $r\geq 3$.
Let $0<a<b<r$ be integers.
We consider the following rational functions
$$
f(u):=\frac{\theta(\pi_0^au)^{r}}{\theta(u)^{r-a}\theta(qu)^a}
=(-u)^a\left(\frac{\theta(\pi_0^au)}{\theta(u)}
\right)^{r}
$$
and
$$
g(u):=\frac{\theta(\pi_0^bu)^{r}}{\theta(u)^{r-b}\theta(qu)^b}
=(-u)^b\left(\frac{\theta(\pi_0^b u)}{\theta(u)}
\right)^{r}
$$
on $E_L:=E_K\ot_K L$.
It is easy to see that the symbol
\begin{equation}\label{ts0}
\xi_L:=\left\{
\frac{f(u)}{f(\pi_0^{-b})},
\frac{g(u)}{g(\pi_0^{-a})}
\right\}
\end{equation}
is contained in the $K$-cohomology group $H^0_\Zar(E_L,{\cal K}_2)$.

\begin{prop}\label{sa}
Put
\begin{equation}\label{notan}
S(\alpha)\os{{\mathrm def}}{=}\prod_{k=1}^\infty
\left(\frac{1-\alpha q^k}{1-\alpha^{-1} q^k}\right)^k
\quad (\alpha\in \ol{K}^*-q^{\Z}).
\end{equation}
Then
$$
\tau^\et_\infty(\xi_L)=
(-1)^{a(r-b)}
\pi_0^{a(b-a)(b-r)}
\left(
\frac{\theta(\pi_0^{b})^{b}}
{\theta(\pi_0^{b-a})^{b-a}\theta(\pi_0^{a})^{a}}
\right)^r
\left(
\frac{S(\pi_0^{b})}{S(\pi_0^{b-a})S(\pi_0^{a})}
\right)^{r^2}
\in L^*/n.
$$
\end{prop}
\begin{pf}
We denote by $\hat{\xi}_L\in K_2(L\langle u \rangle)/n$ the image
of the symbol $\xi_L$. Due to Theorem \ref{formula}, we have
$\tau^\et_\infty(\xi_L)= \hat{\tau}^\et_\infty(\hat{\xi}_L)$. Note
$$
\prod_{k>N}(1-\pi_0^aq^ku)(1-\pi_0^{-a}q^ku^{-1})\in
(K\langle u \rangle^{\ast})^{n}
$$
for sufficiently large $N\gg \nu $.
Therefore we see
\begin{equation}\label{mod1}
\hat{\xi}_L\equiv\left\{
\frac{(-u)^a}{f(\pi_0^{-b})}
\left(\frac{\theta_N(\pi_0^au)}{\theta_N(u)}
\right)^{r},
\frac{(-u)^b}{g(\pi_0^{-a})}
\left(\frac{\theta_N(\pi_0^bu)}{\theta_N(u)}
\right)^{r}
\right\}
\mod n K_2(L\langle u \rangle)
\end{equation}
where we put
$$
\theta_N(u):=(1-u)\prod_{k=1}^N(1-q^ku)(1-q^ku^{-1}).
$$
The right hand side of \eqref{mod1} comes from $K_2(K(u))$,
so that we can calculate $\hat{\tau}^\et_\infty(\hat{\xi}_L)$ by the
tame symbol (Theorem \ref{formula}).
The following are straightforward:
\begin{equation}\label{formula1}
\hat{\tau}^\et_\infty\left\{
u,c
\right\}
=
c^{-1}
\end{equation}
\begin{equation}\label{formula2}
\hat{\tau}^\et_\infty\left\{
\theta_N(\pi_0^iu),c
\right\}
=
1
\end{equation}
\begin{equation}\label{formula3}
\hat{\tau}^\et_\infty\left\{
\theta_N(\pi_0^iu),u
\right\}
=
1
\end{equation}
\begin{equation}\label{formula4}
\hat{\tau}^\et_\infty\left\{
\theta_N(\pi_0^iu),\theta_N(\pi_0^ju)
\right\}
=
S(\pi_0^{i-j})
\end{equation}
for all $0\leq i,~j<r$.
Using the aboves, we have
$$
\hat{\tau}^\et_\infty
\{f,g\}=(-1)^{ab}
\left(
\frac{S(\pi_0^{b})}{S(\pi_0^{b-a})S(\pi_0^{a})}
\right)^{r^2},
$$
and
\begin{align*}
\tau^\et_\infty(\xi_L)&=
\hat{\tau}^\et_\infty
(\hat{\xi}_L)\\
&=
\frac{f(\pi_0^{-b})^{-b}}{g(\pi_0^{-a})^{-a}}
\cdot \hat{\tau}^\et_\infty
\{f,g\}\\
&=
(-1)^{a(r-b)}
\pi_0^{a(b-a)(b-r)}
\left(
\frac{\theta(\pi_0^{b})^{b}}
{\theta(\pi_0^{b-a})^{b-a}\theta(\pi_0^{a})^{a}}
\right)^r
\left(
\frac{S(\pi_0^{b})}{S(\pi_0^{b-a})S(\pi_0^{a})}
\right)^{r^2} .
\end{align*}
\end{pf}
\begin{cor}\label{bysato}
Let $N_{L/K}:H^0_\Zar(E_L,{\cal K}_2)\to
H^0_\Zar(E_K,{\cal K}_2)$ be the norm map.
Put $\xi=N_{L/K}(\xi_L)$.
Then $o_K(\xi)$ is a nonzero integer.
In particular the cokernel of the map \eqref{sato} is finite.
\end{cor}
\begin{pf}
$\tau^\et_\infty(\xi_L)$ can be written as
$$
\tau^\et_\infty(\xi_L)=
(-1)^{a(r-b)}
\pi_0^{a(b-a)(b-r)}
\left(
\frac{M_b}{M_{b-a}M_a}\right)^r\in L^*/n
$$
where
$$
M_i:=
\theta(\pi_0^{i})^{i}S(\pi_0^i)^r
=\prod_{n=1}^\infty
\frac{(1-\pi_0^{nr-r+i})^{nr-r+i}}{(1-\pi_0^{nr-i})^{nr-i}}.
$$
Therefore we have $o_L(\xi_L)=a(b-a)(b-r)$.
Denote by $f$ the degree of the residue field of $L$
over the residue field of $K$.
Since the diagram
$$
\begin{CD}
H^0_\Zar(E_L,{\cal K}_2) @>{o_L}>> \hat{\Z}\\
@V{N_{L/K}}VV@VV{\text{mult. by } f}V\\
H^0_\Zar(E_K,{\cal K}_2) @>{o_K}>> \hat{\Z}
\end{CD}
$$
is commutative we have $o_K(\xi)=fa(b-a)(b-r)$. This is a nonzero integer.
\end{pf}
\begin{cor}[T. Sato \cite{TSato}]\label{bysator}
The $l$-adic regulator
$H^0_\Zar(E_K,{\cal K}_2)\ot\Q_l\to H^2_\et(E_K,\Q_l(2))$ is surjective
for $l\not=p$.
\end{cor}
\begin{pf}
In the same way as the $p$-adic case,
we can show that
there is an exact
sequence
$$
0\lra H^1_\et(K,\Q_l(2))\lra H^2_\et(E_K,\Q_l(2))\lra
H^1_\et(K,\Q_l(1))\lra 0.
$$
Since $l\not=p$, it follows from the Euler-Poincar\'e characteristic
(cf. \cite{galois} II 5.7.)
that we have $H^1_\et(K,\Q_l(2))=0$ and
$H^1_\et(K,\Q_l(1))\os{\cong}{\to} \Q_l$.
The composition
$$
H^0_\Zar(E_K,{\cal K}_2)\lra H^2_\et(E_K,\Q_l(2))
\os{\cong}{\lra}H^1_\et(K,\Q_l(1))\os{\cong}{\lra} \Q_l
$$
coincides with the composition of $o_K$ and the natural map
$\hat{\Z}\to \Q_l$. Therefore the surjectivity follows from
Corollary \ref{bysato}.
\end{pf}
\begin{rem}[T.Sato's thesis]\label{thesis}
The proof of Corollary \ref{bysator} is different from
his one in \cite{TSato}. His proof is done in the following way.

Let ${\cal E}/R$ be the N\'eron model of $E_L/L$, and
$\G_{m,k}\subset {\cal E}_k$ the identity component of
the special fiber.
Let
$$
\partial=\partial_2\partial_1:K_2(E_L) \os{\partial_1}{\lra} K_1(\G_{m,k})
\os{\partial_2}{\lra} K_0(k)=\Z
$$
be the composition
where $\partial_i$ is the boundary map coming from
the localization sequence of $K$-theory.
(It seems $\partial=o_L$ under the inclusion
$\Z\hookrightarrow\hat{\Z}$.)
In his thesis, T.Sato constructed the symbol $\xi_L$ and showed
$\partial(\xi_L)=a(b-a)(b-r)$ (cf. \cite{scholl} for more calculation
of the boundary).
This implies
$$
{\mathrm corank}~ H^0(E_K,{\cal K}_2)\ot\Q_l/\Z_l \geq 1,
$$
and thus he obtained the surjectivity of the $l$-adic regulator
$K_2(E_K)\ot\Q_l\to H^2_\et(E_K,\Q_l(2))$ and
the finiteness of $l$-power torsion part
$K_1(E_K)[l^\infty]$ for $l\not=p$
by using Suslin's exact sequence (cf. Proposition \ref{equivthm} below).

The map $\partial$ is enough to study
the $l$-adic regulator on $K_2$.
However, it is not enough to study the $p$-adic regulator.
In fact, our map
$\tau^\et_\infty$ plays an essential role in the next step.
\end{rem}

\subsection{Proof of Part I : Step 2}
Put $U_{K}:=(1+\pi_K R)/p^\nu\subset K^*/p^\nu$ and
$$
U_{K}^{\tau^\et_\infty}:=U_{K}\cap\tau^\et_\infty(
H^0(E_K,{\cal K}_2))
$$
Our next step is to show that the cardinality of
$U_{K}/U^{\tau^\et_\infty}_{K}$ has an
upper bound which does not depend on $\nu$.
Step 1 and Step 2 immediately imply (Part I).

It follows from the norm map that
we may replace $K$ with any finite extension $L$ of $K$.
Thus we may assume that there is a uniformizer $\pi_0\in K^*$
such that $\pi_0^r=q$.
The proof is done in the following steps:

\smallskip

{\bf (Step 2-1)} \quad
Let $i~m\geq 1$ be any integers and $\zeta_m$
any $m$-th root of unity. Suppose $\zeta_m\in K^*$.
Then we have
$
(1-\zeta_m\pi_0^i)^{mi}\in U^{\tau^\et_\infty}_{K}.
$

{\bf (Step 2-2)} \quad
Let $V_{K}\subset U_{K}$ be the subgroup generated by
all $(1-\zeta_m\pi_0^i)^{mi}$ where $i\geq 1$ and $\zeta_m\in K^*$
are roots of unity with $\zeta^m_m=1$.
(By Step 2-1, we have $V_{K}\subset U^{\tau^\et_\infty}_{K}\subset U_K$.)
Then the cardinality of $U_{K}/V_{K}$
has an upper bound which does not depend on $\nu$.

\medskip

To prove the above steps, we use the following lemmas.
\begin{lem}\label{f0}
Let $m_1,~m_2\geq 1$ be integers, and
$\zeta_i\in K^*$ $m_i$-th roots of unity with $\zeta_1\not=\zeta_2$.
Put
$$
f(u):=\left(\frac{\theta(\zeta_1^{-1}u)}{\theta(u)}\right)^{m_1},
\quad
g(u):=\left(\frac{\theta(\zeta_2^{-1}u)}{\theta(u)}\right)^{m_2}.
$$
Then
$$
\tau^\et_\infty\left\{
\frac{f(u)}{f(\zeta_2)},
\frac{g(u)}{g(\zeta_1)}
\right\}
=
\left(
\frac{S(\zeta_1^{-1}\zeta_2)}{S(\zeta_1^{-1})S(\zeta_2)}
\right)^{m_1m_2}.
$$
Here $S(\alpha)$ is as in \eqref{notan}.
\end{lem}
\begin{lem}\label{f1}
Let $m\geq 1$ and $1\leq b < a$ be integers, and
$\zeta\in K^*$ a root of unity with $\zeta^m=1$.
Suppose that there is a
$q_0\in K^*$ such that $q_0^a=q$.
Put
$$
f(u):=
\frac{\theta(q_0^{-b}u)^a}{\theta(u)^{a-1}\theta(q^{-b}u)},
\quad
g(u):=\left(\frac{\theta(\zeta^{-1}u)}{\theta(u)}\right)^{m}.
$$
Then
\begin{align*}
\tau^\et_\infty\left\{
\frac{f(u)}{f(\zeta)},
\frac{g(u)}{g(q_0^b)}
\right\}
&=
\left(\frac{S(q_0^{-b}\zeta)}{S(\zeta)S(q_0^{-b})}
\right)^{ma}
\left(\frac{\theta(q_0^{b})}{\theta(q_0^{b}\zeta^{-1})}
\right)^{mb}\\
&=
\left(
S(\zeta)^{-a}
\left(
\frac{1-q_0^b}{1-\zeta^{-1}q_0^b}\right)^b
\prod_{k=1}^\infty
\left(
\frac{1-q_0^{b}q^k}{1-\zeta^{-1}q_0^{b}q^k}
\right)^{ak+b}
\left(
\frac{1-\zeta q_0^{-b}q^k}{1-q_0^{-b}q^k}
\right)^{ak-b}\right)^m.
\end{align*}
\end{lem}
The proofs of Lemmas \ref{f0} and \ref{f1} are similar to the one
of Proposition \ref{sa}.

\subsection{Proof of Step 2-1.}
\begin{lem}\label{step2-1}
Let $\zeta\in K^*$ be any root of unity with
$\zeta^m=1$. Then
$S(\zeta)^m\in  U^{\tau^\et_\infty}_{K}$.
\end{lem}
\begin{pf}
Let $\mu\in \ol{K}^*$ be a primitive $N$-th root of unity
with $(N,pm)=1$. Put $L=K(\mu)$ and $E_L=E_K\ot L$.
By Lemma \ref{f0}, we have
\begin{equation}\label{step2-11}
\left(\frac{S(\zeta\mu)}{S(\zeta)S(\mu)}\right)^{mN}
\in  U^{\tau^\et_\infty}_{L}
~(\subset L^*/p^\nu).
\end{equation}
Since $p\not\vert N$, we have
\begin{equation}\label{step2-12}
\left(\frac{S(\zeta\mu)}{S(\zeta)S(\mu)}\right)^m
\in  U^{\tau^\et_\infty}_{L}
~(\subset L^*/p^\nu).
\end{equation}
Applying the norm map for $L/K$, we have
\begin{equation}\label{step2-13}
\left(\prod_{i=0}^{N-1}\frac{S(\zeta\mu^i)}{S(\zeta)S(\mu^i)}
\right)^m
=\left(S(\zeta)^{-N}
\prod_{i=0}^{N-1}\frac{S(\zeta\mu^i)}{S(\mu^i)}\right)^m
\in  U^{\tau^\et_\infty}_{K}
~(\subset K^*/p^\nu).
\end{equation}
Choosing a sufficiently large $N\gg 1$ with $(N,pm)=1$, we have
$$
\prod_{i=0}^{N-1}S(\zeta\mu^i)=
\prod_{k=1}^{\infty}
\left(\frac{1-\zeta^Nq^{Nk}}{1-\zeta^{-N}q^{Nk}}\right)^k
\equiv 1
\mod (K^*)^{p^\nu}
$$
and $\prod_{i=0}^{N-1}S(\mu^i)\equiv 1$. Thus we have
\begin{equation}\label{step2-14}
S(\zeta)^{-Nm}
\in  U^{\tau^\et_\infty}_{K}.
\end{equation}
Since $p\not\vert N$, we have $S(\zeta)^{m}
\in  U^{\tau^\et_\infty}_{K}$.
\end{pf}
\begin{lem}\label{step2-2}
$(1-\pi^i_0)^i\in  U^{\tau^\et_\infty}_{K}$
for all $i\geq 1$.
\end{lem}
\begin{pf}
Let $m,~a\geq 1$ be integers.
Let $\zeta\in\ol{K}^*$ be a primitive $m$-th root of unity.
Take $q_0\in\ol{K}^*$ such that $q_0^a=\pi_0$ (and hence $q_0^{ar}=q$).
We put $L_1=K(\zeta)\subset L_2=K(q_0,\zeta)$ and $E_{L_i}
=E_K\ot L_i$.
Due to Lemmas \ref{f1} and \ref{step2-1},
we have
\begin{equation}\label{step2-21}
\left(
\left(\frac{1-q_0^b}{1-\zeta^{-1}q_0^b}\right)^b
\prod_{k=1}^\infty
\left(
\frac{1-q_0^{b}q^k}{1-\zeta^{-1}q_0^{b}q^k}
\right)^{ark+b}
\left(
\frac{1-\zeta q_0^{-b}q^k}{1-q_0^{-b}q^k}
\right)^{ark-b}
\right)^m
\in  U^{\tau^\et_\infty}_{L_2}
\end{equation}
for any $1\leq b< ar$.
Suppose that $a$ is large enough and
$(a,bm)=1$.
Taking the norm map for $L_2/L_1$, we have
\begin{equation}\label{step2-22}
\left(
\left(\frac{1-\pi_0^b}{1-\zeta^{-a}\pi_0^b}\right)^b
\prod_{k=1}^\infty
\left(
\frac{1-\pi_0^{b}q^{ak}}{1-\zeta^{-a}\pi_0^{b}q^{ak}}
\right)^{ark+b}
\left(
\frac{1-\zeta^a \pi_0^{-b}q^{ak}}{1-\pi_0^{-b}q^{ak}}
\right)^{ark-b}
\right)^m
\in  U^{\tau^\et_\infty}_{L_1}.
\end{equation}
Since $a\gg 1$, we have
\begin{equation}\label{step2-23}
\left(\frac{1-\pi_0^b}{1-\zeta^{-a}\pi_0^b}\right)^{mb}
\in  U^{\tau^\et_\infty}_{L_1}.
\end{equation}
Suppose further $m\gg 1$ and $p\not\vert m$.
Then by taking the norm map for $L_1/K$,
we have
\begin{equation}\label{step2-24}
\left(\frac{(1-\pi_0^b)^m}{1-\pi_0^{mb}}\right)^{mb}
\equiv
(1-\pi_0^b)^{m^2b}
\in U^{\tau^\et_\infty}_{K}.
\end{equation}
Since $p\not\vert m$, we have
$(1-\pi_0^b)^{b}
\in U^{\tau^\et_\infty}_{K}$.
\end{pf}
\begin{lem}\label{step2-2cor}
$(1-\zeta\pi^i_0)^{mi}\in  U^{\tau^\et_\infty}_{K}$
for all $i\geq 1$ and all roots of unity $\zeta\in K^*$
such that $\zeta^m=1$.
\end{lem}
\begin{pf}
In the proof of Lemma \ref{step2-2},
the same argument works until \eqref{step2-23}.
Thus we have
\begin{equation}\label{prestep2-3}
\left(\frac{1-\pi_0^b}{1-\zeta\pi_0^b}\right)^{mb}
\in U^{\tau^\et_\infty}_{K}
\end{equation}
for all $b\geq 1$.
By Lemma \ref{step2-2}, we have
$(1-\zeta\pi_0^{b})^{mb}
\in U^{\tau^\et_\infty}_{K}$.
\end{pf}
Lemma \ref{step2-2cor} completes the proof of Step 2-1.

\subsection{Proof of Step 2-2.}
Let
$U^i_K$ be the subgroup of $K^*/p^\nu$ generated by
$1+\pi_K^iR$ and we put $V^i_K:=V_K\cap U^i_K$.
By definition $U^1_K=U_K$, and $U^i_K=0$ for $i\gg \nu$.
Let $e$ be the ramified index of $K/\Q_p$ (i.e. $\pi_K^eR=pR$).
Then we show that the map
\begin{equation}\label{2-21}
V^i_K/V_K^{i+1} \lra U_K^i/U_K^{i+1}
\end{equation}
is surjective for all
$i\geq e^2+2e$.
This implies $U^i_K=V^i_K$ for
$i\geq e^2+2e$, and hence we obtain an upper bound of $U_K/V_K$
which does not depend on $\nu$.

Write $i=ek+l$ with $e+1\leq l\leq 2e$.
Since $i\geq e^2+2e$, we have $k\geq e$.
Let $a\in R^*$ be any invertible element.
Since $i/e>k+1/(p-1)$, there is an invertible element $a'\in R^*$
such that $1+a\pi_0^i=(1+a'\pi_0^l)^{p^k}$.
It follows from $l\leq 2e \leq p^e< p^{k+1}$ that $\ord_p(l)\leq k$
and therefore we have $1+a\pi_0^i\in V^i_K\cdot (U_K^{l+1})^{p^k}$.
On the other hand, since $l+1\geq e$, we have
$$
\ord_K\left(\binom{p^k}{s}\pi_0^{sl+s}\right)
=e(k-\ord_p(s))+sl+s \geq i+1
$$
for all $1\leq s \leq p^k$. This shows $(U_K^{l+1})^{p^k}\subset
U_K^{i+1}$.
Thus we have $1+a\pi_0^i\in V_K^i\cdot U_K^{i+1}$.



\section{Nodal rational curves and the Bloch groups}\label{partiipre}
Before going to (Part II), we study $K_2$ of nodal rational
curves. To do this, we need $K$-theory and regulator not only for schemes
but also for simplicial schemes.
We work in A. Huber's theory \cite{huber} which suffices for our purpose.
\subsection{Higher $K$-theory of simplicial schemes}
Let $F$ be a field of characteristic zero. We work over the
category $({\mathrm Sch}/F)$ of separated schemes of finite type
over $\Spec F$. Let $\Delta$ be the category of finite sets
$\{0,\cdots,n\}$ with ordering $\leq $. A {\it simplicial scheme}
is a functor from $\Delta^{\mathrm op}$ to $({\mathrm Sch}/F)$. We
write $X_n=X_\bullet(\{0,\cdots,n\})$ for a simplicial scheme
$X_\bullet$. A scheme $X$ is canonically considered as the
simplicial scheme such that $X_n=X$ for $n\geq 0$.

\medskip

The $K$-groups $K(X_\bullet)$
of a simplicial scheme $X_\bullet$ are defined.
They are functorial and agree with
the usual $K$-theory if $X_\bullet$ is a scheme.
We refer \cite{friedlanderEHSS} or \cite{huber} for the details.
Rather than going into the general theory, we
pick up the results which we will use later.
\begin{thm}[\cite{huber} Proposition 18.1.2]
Let $X_\bullet$ be a simplicial scheme. Then there is a natural
spectral sequence
\begin{equation}\label{simpspectral}
E_1^{pq}=
\begin{cases}
0& p<0\\
K_q(X_p)\cap \ker s^0\cap\cdots \ker s^{p-1}
& \text{others}
\end{cases}
\Longrightarrow K_{q-p}(X_\bullet)
\end{equation}
where $s^i:K_q(X_p) \to K_q(X_{p-1})$ are the degeneracy maps.
\end{thm}
A simplicial scheme $X_\bullet$ is called {\it split}
if
$$
N(X_n)=X_n-\bigcup_{s}s(X_{n-1})
$$
is an open and closed subscheme of $X_n$. Here $s:X_{n-1}\to X_n$
runs over all degeneracy maps.
We mostly work over split simplicial schemes with finite combinatorial
dimension, namely simplicial schemes which are split and
such that $N(X_n)$ is empty
for large $n$.
If $X_\bullet$ is a split simplicial scheme with finite combinatorial
dimension, then the spectral sequence \eqref{simpspectral} converges
and $E_1^{pq}=K_q(N(X_p))$.
\begin{thm}[loc.cit. 18]
There is the regulator map
$$
c_{i,j}:K_i(X_\bullet)\lra H^{2j-i}_\et(X_\bullet,\Z/p^\nu(j)),
\quad i,~j\geq 0.
$$
If $F=\C$ and each $X_n$ is nonsingular, then
we also have the regulator map
$$
c_{i,j}^D:K_i(X_\bullet)\lra H^{2j-i}_D(X_\bullet,\Z(j)),
\quad i,~j\geq 0
$$
to the Deligne-Beilinson cohomology. They are functorial and agree
with the usual regulator maps \eqref{chern} when $X_\bullet$ is a
scheme.
\end{thm}

\subsection{Nodal rational curves}\label{nodalsec}
Let $C$ be an irreducible rational curve over $F$
with one node.
It is obtained by attaching $0$ to $\infty$ of $\P^1$.
We denote by $*$ the node of $C$.
Let $f:\P^1 \to C$ be the normalization such that $f^{-1}(*)=
\{0,\infty\}$.
We have the simplicial scheme $C_\bullet$ from $C$ in the usual way,
namely,
$C_0=\P^1$, $C_1=\P^1\coprod \{*\},
\cdots$, and
$d_i:C_1\to C_0$ is defined as the identity on $\P^1$ and
on $\{*\}$,
$d_0=i_0$ the inclusion into $0$ and
$d_1=i_\infty$ the inclusion into $\infty$, etc.
The natural map $C_\bullet\to C$ is a proper hypercovering
(\cite{HodgeIII} 5.3.5 V).
Since $C_\bullet$ is a
split simplicial scheme with finite combinatorial dimension, we have
an exact sequence
\begin{equation}\label{simpspec0}
\cdots \lra K_{i+1}(F)\lra K_i(C_\bullet)\lra K_i(\P^1)
\os{i_0^*-i_\infty^*}{\lra}
K_i(F) \lra \cdots
\end{equation}
from \eqref{simpspectral}. The composition of the natural maps
$K_i(C)\to K_i(C_\bullet)$ and $K_i(C_\bullet)\to K_i(\P^1)$ is
equal to the pull-back $f^*$. Moreover, we claim
$i_0^*-i_\infty^*=0$. In fact, $K_i(\P^1)$ is isomorphic to
$K_i(F)\ot K_0(\P^1) \cong K_i(F)^{\op 2}$. In the commutative
diagram
$$
\begin{CD}
K_i(F)\ot K_0(\P^1)
@>{{\mathrm id}\ot(i^*_0-i^*_\infty)}>>
K_i(F)\ot K_0(F)\\
@V{\cong}VV@VV{\cong}V\\
K_i(\P^1) @>{i^*_0-i^*_\infty}>> K_i(F)
\end{CD}
$$
the above map is clearly zero.
Thus we have $i^*_0-i^*_\infty=0$.
Now the exact sequence \eqref{simpspec0} becomes
\begin{equation}\label{simpnode1}
0\lra K_{i+1}(F)\lra K_i(C_\bullet)\lra K_i(\P^1) \lra 0.
\end{equation}
Put $K_i(C)_0:=\ker (f^*:K_i(C)\to K_i(\P^1))$.
From \eqref{simpnode1}, we have
a natural map
\begin{equation}\label{natural}
\delta:K_i(C)_0 \lra K_{i+1}(F).
\end{equation}
The similar argument also works on \'etale cohomology, and they are
compatible under the regulator maps.
Therefore we have a commutative diagram
\begin{equation}\label{natural2}
\begin{CD}
K_i(C)_0 @>{\delta}>> K_{i+1}(F)\\
@V{c_{i,j}}VV@VV{c_{i+1,j}}V\\
H_\et^{2j-i}(C,\Z/p^\nu(j))_0 @>{\delta_\et}>>
H_\et^{2j-i-1}(F,\Z/p^\nu(j))
\end{CD}
\end{equation}
where we put
$H_\et^{2j-i}(C,\Z/p^\nu(j))_0:=\ker(f^*:
H_\et^{2j-i}(C,\Z/p^\nu(j))
\to
H_\et^{2j-i}(\P^1,\Z/p^\nu(j)))$.

Of particular interest to us is the case $i=j=2$.
Write $C_{\ol{F}}=C\ot_F{\ol{F}}$.
We can easily see
$H_\et^{2}(C,\Z/p^\nu(2))_0=H^1(F,
H_\et^1(C_{\ol{F}},\Z/p^\nu(2)))$ and
$\delta_\et$ is the map defined from the natural isomorphism
$H_\et^1(C_{\ol{F}},\Z/p^\nu(2))\cong
\Z/p^\nu(2)$ up to sign.
As a result, we obtain
\begin{prop}
The diagram
\begin{equation}\label{natural3}
\begin{CD}
K_2(C)_0 @>{\delta}>> K_3(F)\\
@V{\rho_C}VV@VV{c_{3,2}}V\\
H_\et^1(F,H_\et^1(C_{\ol{F}},\Z/p^\nu(2)) @>{\cong}>>
H_\et^1(F,\Z/p^\nu(2))
\end{CD}
\end{equation}
is commutative up to sign. Here $\rho_C$ is given in
\eqref{rho}, and
the isomorphism below is induced
from the natural isomorphism
$H_\et^1(C_{\ol{F}},\Z/p^\nu(2)))\cong
\Z/p^\nu(2)$ up to sign.
\end{prop}
\begin{rem}
To remove the sign ambiguity, we need a careful looking at
the relation between the map $K_3(F)\to K_2(C_\bullet)$
and the isomorphism
$H_\et^1(C_{\ol{F}},\Z/p^\nu(2)))\cong
\Z/p^\nu(2)$. Since it is nothing important for our purpose,
we omit it.
\end{rem}
\subsection{Local ring of the Nodal curve}
Let $\O_{0=\infty}$ be the local ring of $C$ at $*$.
More explicitly, it is given as follows:
\begin{align*}
\O_{0=\infty}&=\left\{f(t)\in F(t)~\vert~
f(0)=f(\infty)\not=\infty
\right\}\\
&=\left\{
\alpha+\beta\frac{t^n+a_1t^{n-1}+\cdots+a_n}
{t^n+b_1t^{n-1}+\cdots+b_n}~\vert~
\alpha,~\beta,~a_i,~b_j\in F,~a_n=b_n\not=0
\right\}.
\end{align*}
Quillen's localization theorem (\cite{Q2},
\cite{Sr} Thm.(9-1)) yields the exact sequence
\begin{equation}\label{quillen1}
0 \lra K_2(C)_{\Q} \lra K_2(\O_{0=\infty})_\Q \os{\tau}{\lra}
\bigoplus_{x\in C-\{*\}}
\kappa(x)^*_\Q
\end{equation}
where $\tau$ is the tame symbol \eqref{tame}.
Moreover, by a theorem of van der Kallen \cite{vdK},
$K_2(\O_{0=\infty})$ is isomorphic to Milnor's $K_2^M(\O_{0=\infty})$.
Thus we can think of $K_2(C)_{\Q}$ being a subgroup of
$K_2^M(\O_{0=\infty})_\Q$.

Let $\O_{0,\infty}$ be the semi-local ring of $\P^1$ at $0$ and $\infty$.
Let $\O_\bullet$ be the simplicial scheme
associated to $\Spec\O_{0=\infty}$.
Similarly to \eqref{simpspec0}, we have
\begin{equation}\label{simpspec1}
\begin{CD}
@.@. K_i^M(\O_{0=\infty})\\
@.@.@VVV\\
\cdots @>>> K_{i+1}(F)@>>> K_i(\O_\bullet)@>>> K_i(\O_{0,\infty})
@>{i_0^*-i_\infty^*}>>
K_i(F) @>>> \cdots
\end{CD}
\end{equation}
where $i_0:\{*\}\to \Spec\O_{0,\infty}$ and $i_\infty:\{*\}\to
\Spec\O_{0,\infty}$ are the inclusions into $0$ and $\infty$ respectively.
\begin{lem}\label{6}
Define the indecomposable $K_3$-group $K_3^{\mathrm ind}(F)$
as the cokernel of the natural map
$K^M_3(F)\to K_3(F)$.
Then the cokernel of $i_0^*-i^*_\infty:
K_3(\O_{0,\infty})\to K_3(F)$ is isomorphic to
$K_3^{\mathrm ind}(F)$.
\end{lem}
\begin{pf}
By Quillen's localization theorem, we have
$$
K_3(\P^1-\{1\}) \lra
K_3(\O_{0,\infty}) \lra \bigoplus_{x\not=0,1,\infty} K_2(\kappa(x)).
$$
Note that $K_3(\P^1-\{1\})=K_3(F)$.
The composition of the maps $K^M_3(\O_{0,\infty})
\to K_3(\O_{0,\infty}) \to \bigoplus_{x\not=0,1,\infty} K_2(\kappa(x))$
is the tame symbol. A direct calculation yields that it is surjective.
This shows that the
map $K^M_3(\O_{0,\infty})\to
K_3(\O_{0,\infty})/K_3(F)$ is surjective.
Therefore the cokernel of $i_0^*-i^*_\infty:
K_3(\O_{0,\infty})\to K_3(F)$ is equal to the cokernel of
$i_0^*-i^*_\infty:
K_3^M(\O_{0,\infty})\to K_3(F)$. It is easy to see that
the image of $i_0^*-i^*_\infty:
K_3^M(\O_{0,\infty})\to K_3(F)$ is the image of $K^M_3(F)$.
This completes the proof.
\end{pf}
We put $K_2^M(\O_{0=\infty})_0=\ker(f^*:
K_2^M(\O_{0=\infty})\to K_2(\O_{0,\infty}))$.
By Lemma \ref{6} and \eqref{simpspec1}, we have
a map
\begin{equation}\label{naturalmi}
K_2^M(\O_{0=\infty})_0\lra K_3^{\mathrm ind}(F).
\end{equation}
It is clearly compatible with \eqref{natural}:
\begin{equation}\label{quillen2}
\begin{CD}
K_2^M(\O_{0=\infty})_0@>>> K^{\mathrm ind}_3(F)\\
@AAA@AAA\\
K_2(C)_0@>{\delta}>> K_3(F).
\end{CD}
\end{equation}
Note that the natural map $K_2(C)_0\ot\Q
\to K_2^M(\O_{0=\infty})_0\ot\Q$
is bijective due to \eqref{quillen1}
and the injectivity of $K_2(\P^1)_\Q\to K_2(\O_{0,\infty})_\Q$.

\subsection{Bloch groups.}
Let $D(F)$ be the free abelian group with basis $[x]$
($x\in F^*-\{1\}$), and $P(F)$ the quotient group of $D(F)$
by the subgroup generated by the following
\begin{equation}\label{scissor}
[x]-[y]+[y/x]-[\frac{1-x^{-1}}{1-y^{-1}}]+[\frac{1-x}{1-y}]
\quad (x\not=y\in F^*-\{1\}).
\end{equation}
The relation \eqref{scissor} is called the
{\it scissors congruence relations}.
Then one can easily derive
the following basic relations in $P(F)\ot\Q$
(cf. \cite{scissor2} \S 5):
\begin{equation}
[x]+[x^{-1}]=0 \quad (x\in F^*),
\end{equation}
\begin{equation}
[x]+[1-x]=0\quad (x\in F^*-\{1\}).
\end{equation}
If $F$ contains a primitive $m$-th root $\zeta$ of unity, then
\begin{equation}
\label{zetarel}
m\sum_{i=1}^m[\zeta^i x]=[x^m].
\end{equation}
A homomorphism
\begin{equation}\label{lambda1}
\lambda:P(F)\lra F^*\wedge F^*, \quad  [x]\mapsto x\wedge(1-x).
\end{equation}
is well-defined.
The kernel of $\lambda$ is called
the {\it Bloch group} which we denote by $B(F)$:
$$
0\lra B(F) \lra P(F) \os{\lambda}{\lra}
 F^*\wedge F^*\lra K_2^M(F) \lra 0.
$$

Using some ideas of Bloch, Suslin proved the following remarkable
theorem.
\begin{thm}[Suslin \cite{suslinicm}]\label{bloch1}
$K_3^{\mathrm ind}(F)_\Q\cong B(F)_\Q$.
\end{thm}
See also related works by Dupont and Sah \cite{scissor2}.
Hereafter we identify $K_3^{\mathrm ind}(F)_\Q$ with $B(F)_\Q$
by the above theorem.
\subsection{Explicit Description of $\delta$}
Passing to the projective limit and tensoring with $\Q$,
we have from \eqref{natural3}
\begin{equation}\label{natural4}
\begin{CD}
K_2(C)_{0,\Q} @>{\delta}>> K_3(F)_{\Q}\\
@V{\rho_C}VV@VV{c_{3,2}}V\\
H_\et^1(F,H_\et^1(C_{\ol{F}},\Q_p(2)) @>{\cong}>>
H_\et^1(F,\Q_p(2)).
\end{CD}
\end{equation}
Note $K_2(C)_{0,\Q}\cong K_2^M(\O_{0=\infty})_{0,\Q}$.
As is well-known, the regulator map $c_{3,2}$ factors through
$K_3^{\mathrm ind}(F)$.
Moreover, $K_3^{\mathrm ind}(F)$ is isomorphic to the Bloch group
$B(F)_\Q$ by Theorem \ref{bloch1}.
We thus have a diagram
\begin{equation}\label{natural4bar}
\begin{CD}
K_2^M(\O_{0=\infty})_{0,\Q} @>{\bar{\delta}}>> B(F)_{\Q}\\
@V{\rho'_C}VV@VV{c_{3,2}}V\\
H_\et^1(F,H_\et^1(C_{\ol{F}},\Q_p(2)) @>{\cong}>>
H_\et^1(F,\Q_p(2))
\end{CD}
\end{equation}
which is commutative up to sign.

We want to describe the map $\bar{\delta}$
explicitly. Unfortunately, it is done only when $F\subset\ol{\Q}$,
because we use Borel's theorem in the proof.
\begin{prop}\label{bloch2}
Suppose $F\subset \ol{\Q}$.
Put $[a,b]:=[a^{-1}b]-[a^{-1}]-[b]\in P(F)$.
Let
$$
\xi=\sum\left\{c\prod_i\frac{1-a_i^{-1}t}{1-b_i^{-1}t},
c'\prod_j\frac{1-c_j^{-1}t}{1-d_j^{-1}t}
\right\}\in K^M_2(\O_{0=\infty})
$$
be a symbol with $\prod a_i/b_i=\prod c_j/d_j=1$.
Assume $\xi\in K^M_2(\O_{0=\infty})_0$.
Then we have
$$
\bar{\delta}(\xi)=
\pm\sum\sum_{i,j} [a_i,c_j]-[a_i,d_j]-[b_i,c_j]+[b_i,d_j]
\in B(F)_\Q.
$$
\end{prop}
\begin{rem}
I believe that the above formula holds
without the assumption ``$F\subset \ol{\Q}$".
\end{rem}
However we use Proposition \ref{bloch2} only for the following special case
(see \S \ref{step4sec} Step 4).
\begin{cor}\label{bloch2cor}
Let $F$ be an arbitrary field of characteristic zero.
Suppose that there are distinct roots of unity $\zeta_1,~\zeta_2\in F$
such that $\zeta_1^{m_1}=\zeta_2^{m_2}=1$.
Let$$
\eta_0:=
\left\{
\left(\frac{1-\zeta_1^{-1}t}{1-t}\right)^{m_1}
\left(\frac{1-\zeta_1^{-1}\zeta_2}{1-\zeta_2}\right)^{-m_1},
\left(\frac{1-\zeta_2^{-1}t}{1-t}\right)^{m_2}
\left(\frac{1-\zeta_2^{-1}\zeta_1}{1-\zeta_1}\right)^{-m_2}
\right\}
$$
be a symbol in $K_2^M(\O_{0=\infty})$.
Then $\eta_0$ is contained in $K_2^M(\O_{0=\infty})_0$,
and
$$
\bar{\delta}(\eta_0)=\pm m_1m_2([\zeta_1\zeta_2^{-1}]-[\zeta_1]
-[\zeta_2^{-1}])\in B(F)_\Q.
$$
\end{cor}
\begin{pf}
Since everything are defined over $\Q(\zeta_1,\zeta_2)$, we may assume
$F=\Q(\zeta_1,\zeta_2)$. Thus we can apply Proposition \ref{bloch2}
if we show $\eta_0\in K_2^M(\O_{0=\infty})_{0}$.

Letting $\eta_0'=f^*\eta_0\in K_2(\O_{0,\infty})$,
we want to show $\eta_0'=0$.
Recall the localization exact sequence
$$
K_2(F)=K_2(\P^1-\{1\}) \lra K_2(\O_{0,\infty}) \lra
\bigoplus_{x\not=0,1,\infty}
\kappa(x)^*.
$$
The composition of the maps
$K_2^M(\O_{0=\infty})\to K_2(\O_{0,\infty}) \to
\bigoplus_{x\not=0,1,\infty}\kappa(x)^*$ is the tame symbol, and
a direct calculation yields the tame image of $\eta_0$ is zero.
Therefore $\eta'_0$ is in the image of $K_2(F)$. We have
\begin{equation}\label{numbretor}
\eta'_0=\eta'_0\vert_{t=0}=
\left\{
\left(\frac{1-\zeta_1^{-1}\zeta_2}{1-\zeta_2}\right)^{-m_1},
\left(\frac{1-\zeta_2^{-1}\zeta_1}{1-\zeta_1}\right)^{-m_2}
\right\}
\end{equation}
in $K_2(\O_{0,\infty})$.
We can see that the right hand side of \eqref{numbretor} is zero
in the following way.
\begin{align*}
\text{R.H.S of \eqref{numbretor}}&=
m_1m_2
\left\{
\frac{1-\zeta_1^{-1}\zeta_2}{1-\zeta_2},
\frac{1-\zeta_2^{-1}\zeta_1}{1-\zeta_1}
\right\}\\
&=m_1m_2
\left\{
\frac{\zeta_1-\zeta_2}{1-\zeta_2},
\frac{\zeta_2-\zeta_1}{1-\zeta_1}
\right\}\\
&=m_1m_2
\left\{1-x, 1-x^{-1}
\right\}\quad
(x:=\frac{1-\zeta_1}{1-\zeta_2})
\\
&=0.
\end{align*}
\end{pf}
\subsection{Proof of Proposition \ref{bloch2}}\label{b2}
We prove the assertion by using
the complex regulators.
For a complex place $\sigma:F\hookrightarrow \C$,
we denote by $c_\sigma$ the composition of
$K_3^{\mathrm ind}(F)\to K_3^{\mathrm ind}(\C)$ and the complex
regulator $c^D_{3,2}:K_3^{\mathrm ind}(\C)\to \R$.
Note $K_3^{\mathrm ind}(F)_\Q=K_3(F)_\Q$ for any number field $F$.
Borel's theorem asserts the isomorphism
$$
K_3^{\mathrm ind}(F)\ot_\Z\R \os{\cong}{\lra} \R^{r_2},\quad
x\longmapsto (\cdots,c_\sigma(x),\cdots).
$$
Put $C_\sigma=C\ot_{F,\sigma}\C$,
$\O_{0=\infty,\sigma}=\O_{0=\infty}\ot_{F,\sigma}\C$ etc.
We denote by ${\rho'_\sigma}$ the composition of
$K^M_2(\O_{0=\infty})_0\to K^M_2(\O_{0=\infty,\sigma})_0$
and the complex regulator map $
K^M_2(\O_{0=\infty,\sigma})_0 \to
\Ext_{\mathrm MHS}(\R,H^1(C_{\sigma},\R(2))$.
Similarly to \eqref{natural4bar}, we have a diagram
\begin{equation}\label{natural3c}
\begin{CD}
K^M_2(\O_{0=\infty})_0 @>{\bar{\delta}}>> B(F)\\
@V{\rho'_\sigma}VV@VV{c_\sigma}V\\
\Ext_{\mathrm MHS}(\R,H^1(C_{\sigma},\R(2)) @>{i}>{\cong}>
\R
\end{CD}
\end{equation}
which is commutative up to sign.
Here the isomorphism $i$ is induced from
the isomorphism $H^1(C_{\sigma},\Z(2))\cong \Z(2)$.
Due to Borel's theorem, it suffices to show
\begin{equation}\label{borel1}
i\rho'_\sigma(\xi)=\pm c_\sigma\bar{\delta}(\xi)=
\pm\sum\sum_{i,j} c_\sigma[a_i,c_j]-c_\sigma[a_i,d_j]
-c_\sigma[b_i,c_j]+c_\sigma[b_i,d_j]\in \R
\end{equation}
for all complex places $\sigma$.

The map $c_{\sigma}$ in \eqref{natural3c}
is given by the {\it Bloch-Wigner function}
$D_2$ (\cite{irvine}):
\begin{equation}\label{borel2}
c_{\sigma}[x]= D_2(\sigma(x)), \quad x\in F^*-\{1\}.
\end{equation}
Here $D_2$ is defined in the following way.
$$
D_2(x)={\mathrm arg}(1-x)\log\vert x\vert-
{\mathrm Im}\int^x_0\log(1-t)\frac{dt}{t}.
$$
This is a singled valued function on $\C-\{0,1\}$.
On the other hand,
the map $\rho'_{\sigma}$
in \eqref{natural3c} is given in the following way.
Let $\sum \{f,g\}$ be a symbol
which is contained in
$K_2^M(\O_{0=\infty})_{0}$.
We denote by $f^\sigma$ the image of $f$ in
$K_2^M(\O_{0=\infty,\sigma})_{0}$.
Choose a path $\gamma\subset \P^1(\C)$ from $0$ to $\infty$
which does not meet either poles or zeros of $f^\sigma$ and
$g^\sigma$, and
such that its homotopy class $[\gamma]$
is a generator of $\pi_1(C_\sigma,*)$.
Then $\rho_{\sigma}$ is given by
\begin{equation}\label{borel3}
i\rho_{\sigma}
(\sum \{f,g\})
=\sum\int_\gamma\log\vert f^\sigma\vert
d{\mathrm arg}(g^\sigma)
-\log\vert g^\sigma\vert d{\mathrm arg}(f^\sigma).
\end{equation}
One can easily check that $\rho_{\sigma}$
does not depends on the choice
of $\gamma$.

Now a direct calculation using \eqref{borel2} and \eqref{borel3}
yields
\eqref{borel1}.
Left to the reader for the details.


\section{Proof of Theorem \ref{mainthm} : Part II}\label{partiipf}
In this section, we prove
\begin{description}
\item[(Part II)]
$\hsymbol
\supset mH_{\mathrm ab}$
for some $m\not=0$ which does not depend on $\nu$
(See \eqref{mot} and \eqref{ab} for the notations.)
\end{description}
\subsection{Proof of Part II : Step 1}
We consider the Tate curve
$E_{n,K}=K^*/q^{n\Z}$ with the period $q^n$ for an
integer $n\geq 1$.
Recall the diagram (cf. \S \ref{rwsec}):
\begin{equation}\label{rw}
\begin{CD}
@. H^0_\Zar(E_{n,K},{\cal K}_2)/p^\nu@>{=}>>
H^0_\Zar(E_{n,K},{\cal K}_2)/p^\nu\\
@.@V{\rho}VV@VV{\tau^\et_\infty}V\\
H^1_\et(K,\Z/p^\nu(2))
@>{a}>>
H^1_\et(K,H^1_\et(E_{n,\ol{K}},\Z/p^\nu(2)))
@>{b}>>
K^*/p^\nu.
\end{CD}
\end{equation}
By (Part I), the cardinality of the cokernel of
$\tau^\et_\infty$ has an upper bound which
does not depend on $\nu$.
The kernel of $a$ is dominated by $H^0(K,\Z/p^\nu(1))$
whose order is at most $N:=\sharp K^*[p^\infty]$.
Let $m_i\geq 1$ ($i=1,~2$) be integers, and
$\zeta_i\in K^*$ $m_i$-th roots of unity with $\zeta_1\not=\zeta_2$.
Let
$$
f(v):=\left(\frac
{\theta(\zeta^{-1}_1v,q^n)}{\theta(v,q^n)}\right)^{m_1},
\quad
g(v):=\left(\frac
{\theta(\zeta_2^{-1}v,q^n)}{\theta(v,q^n)}\right)^{m_2}
$$
be rational functions on $E_{n,K}$, where $\theta(v,q^n)$
is the theta function with the period $q^n$.
Then we consider a symbol
$$
\eta:=\left\{
\frac{f(v)}{f(\zeta_2)},
\frac{g(v)}{g(\zeta_1)}
\right\}
\in H^0_\Zar(E_{n,K},{\cal K}_2)/p^\nu.
$$
By Lemma \ref{f0},
we have $\tau^\et_\infty(\eta)=1$ when
${\mathrm ord}_pq^n>\nu+1/(p-1)$.
Thus we get a class
$$
\wt{\rho(\eta)}\in H^1_\et(K,\Z/p^\nu(2))
$$
such that $a(\wt{\rho(\eta)})=\rho(\eta)$.

On the other hand, let
$C:=\P^1_K/0\sim \infty$ the nodal curve over $K$
which is obtained by attaching
the two points $0$ and $\infty$ (cf. \S \ref{nodalsec}).
We put
$$
\eta_0:=
\left\{
\left(\frac{1-\zeta_1^{-1}t}{1-t}\right)^{m_1}
\left(\frac{1-\zeta_1^{-1}\zeta_2}{1-\zeta_2}\right)^{-m_1},
\left(\frac{1-\zeta_2^{-1}t}{1-t}\right)^{m_2}
\left(\frac{1-\zeta_2^{-1}\zeta_1}{1-\zeta_1}\right)^{-m_2}
\right\}
$$
a symbol in $H^0_\Zar(C,{\cal K}_2)/p^\nu$.
Let
$$
\rho_{C}:H^0_\Zar(C,{\cal K}_2)/p^\nu
\lra H^1_\et(K,H^1_\et(C_{\ol{K}},\Z/p^\nu(2)))
\cong H^1_\et(K,\Z/p^\nu(2))$$
be the regulator as in \eqref{rhoz}.
Thus we get a class
$$
\rho_C(\eta_0)
\in H^1_\et(K,\Z/p^\nu(2)).
$$
\begin{thm}\label{partii}
Let $N$ be the cardinality of $K^*[p^\infty]$.
Suppose ${\mathrm ord}_pq^n
\geq 2\nu+3$ if $p\geq 3$ and
${\mathrm ord}_pq^n\geq 2\nu+5$ if $p=2$.
Then we have
$$
N\cdot \wt{\rho(\eta)}=\pm N\cdot \rho_C(\eta_0)
\in H^1_\et(K,\Z/p^\nu(2)).
$$
\end{thm}
Note that $N\cdot \wt{\rho(\eta)}$ does not depend on the choice of
$\wt{\rho(\eta)}$.
\begin{rem}
The above equality seems true only if ${\mathrm ord}_pq^n
> \nu+1/(p-1)$.
\end{rem}

\def\E{{\cal E}}
\subsection{Proof of Theorem \ref{partii}.}
With an indeterminant $s$, we put
$$
\O_i:=R[[(q^is)]]\subset R[[s]] \quad (i\geq 0).
$$
Since $\O_i$ is isomorphic to $R[[t]]$ as ring,
it is a complete local ring whose maximal ideal is $(\pi_K,q^is)$.
Moreover it is a unique factorization domain
(i.e. any ideal of height 1 is a principal ideal).

\smallskip

Let $A_i:=\O_i[q^{-1},s^{-1}]\subset R[[s]][q^{-1},s^{-1}]$.
Note $\O_i[q^{-1},s^{-1}]=\O_i[\pi_K^{-1},(q^is)^{-1}]$ and hence
$A_i\cong R[[s]][
\pi_K^{-1},s^{-1}]$.
Let $\E_i$ be the Tate curve over $A_i$ with the period $q^is$
\begin{equation}\label{defnring}
\pi:\E_i \lra \Spec A_i.
\end{equation}
Since $\pi$ is a projective and smooth morphism,
the regulator map
$$
H^0_\Zar(\E_i,{\cal K}_2)/p^\nu \lra H^2_\et(\E_i,\Z/p^\nu(2))
$$
gives rise to a map
$$
\rho_i:H^0_\Zar(\E_i,{\cal K}_2)/p^\nu \lra H^1_\et(
A_i,R^1\pi_*\Z/p^\nu(2)).
$$
Moreover, we have an exact sequence
\begin{equation}\label{exbef}
0\lra\Z/p^\nu(2)
\lra
R^1\pi_*\Z/p^\nu(2)
\lra
\Z/p^\nu(1)\lra 0
\end{equation}
of \'etale sheaves on $\Spec A_i$
similarly to \eqref{van0}.
Therefore we have
\begin{equation}\label{ex}
H^1_\et(A_i,\Z/p^\nu(2))
\os{\alpha_i}{\lra}
H^1_\et(A_i,R^1\pi_*\Z/p^\nu(2))
\os{\beta_i}{\lra}
H^1_\et(A_i,\Z/p^\nu(1)).
\end{equation}
The kernel of $\alpha_i$ is dominated by $
H^0_\et(A_i,\Z/p^\nu(1))=H^0(K,\Z/p^\nu(1))$,
which is finite of order $\leq N$.
\begin{lem}\label{Hilbert90}
$H^1_\et(A_i,\Z/p^\nu(1))\cong A_i^*/p^\nu
\cong (\O_i^* \times \pi_K^\Z\times s^\Z)/p^\nu$.
\end{lem}
\begin{pf}
Since $R[[s]]$ is a unique factorization domain,
so is $A_i\cong R[[s]][
\pi_K^{-1},s^{-1}]$.
Then the assertion follows from Hilbert 90.
\end{pf}
Let
$$
f_\E(v):=\left(\frac
{\theta(\zeta^{-1}_1v,q^is)}{\theta(v,q^is)}\right)^{m_1},
\quad
g_\E(v):=\left(\frac
{\theta(\zeta_2^{-1}v,q^is)}{\theta(v,q^is)}\right)^{m_2}
$$
be rational functions on $\E_i$.
Consider the symbol
$$
\eta_\E:=\left\{
\frac{f_\E(v)}{f_\E(\zeta_2)},
\frac{g_\E(v)}{g_\E(\zeta_1)}
\right\}
\in H^0(\E_i,{\cal K}_2)/p^\nu.
$$
For an integer $m\geq 1$
we put by $s_m:A_i\to K$ the $R$-ring homomorphism
given by $s\mapsto q^m$.
Then we have
$$
\E_i\ot_{s_m}K=E_{i+m,K}=K^*/(q^{i+m})^\Z,
$$
$$
s_m^*f_\E(v)=\left(\frac
{\theta(\zeta^{-1}_1v,q^{i+m})}{\theta(v,q^{i+m})}\right)^{m_1},
\quad
s_m^*g_\E(v)=\left(\frac
{\theta(\zeta_2^{-1}v,q^{i+m})}{\theta(v,q^{i+m})}\right)^{m_2}.
$$
\begin{lem}\label{left}
Put
$$
S_{\E_i}(\alpha)\os{{\mathrm def}}{=}\prod_{k=1}^\infty
\left(\frac{1-\alpha (q^is)^k}{1-\alpha^{-1} (q^is)^k}\right)^k
\quad (\alpha\in \O_i^*).
$$
Then we have
\begin{equation}\label{f01e}
(\beta_i\rho_i)(\eta_\E)=
\left(
\frac{S_{\E_i}(\zeta_1^{-1}\zeta_2)}{S_{\E_i}(\zeta_1^{-1})S_{\E_i}(\zeta_2)}
\right)^{m_1m_2}\in
A_i^*/p^\nu.
\end{equation}
\end{lem}
\begin{pf}
Put
$$
S_{l}(\alpha)\os{{\mathrm def}}{=}\prod_{k=1}^\infty
\left(\frac{1-\alpha (q^{l})^k}{1-\alpha^{-1} (q^{l})^k}\right)^k,
\quad l\geq1.
$$
By Lemma \ref{f0} we have
\begin{align*}
s_m\left((\beta_i\rho_i)(\eta_\E)\right)&=
\tau_\infty^\et
\left\{
\frac{s_m^*f_\E(v)}{s_m^*f_\E(\zeta_2)},
\frac{s_m^*g_\E(v)}{s_m^*g_\E(\zeta_1)}
\right\}\\
&=
\left(
\frac{S_{i+m}(\zeta_1^{-1}\zeta_2)}{S_{i+m}(\zeta_1^{-1})S_{i+m}(\zeta_2)}
\right)^{m_1m_2}\\
&=s_m\left(
\frac{S_{\E_i}(\zeta_1^{-1}\zeta_2)}{S_{\E_i}(\zeta_1^{-1})S_{\E_i}(\zeta_2)}
\right)^{m_1m_2}
\in
K^*/p^\nu
\end{align*}
for any $\nu\geq 1$ and $m\geq 1$.
This implies \eqref{f01e} because of the injectivity of
$$
\prod_{m\geq 1}s_m:
\plim{\nu}A_i^*/p^\nu \lra
\prod_{m\geq 1}\plim{\nu}K^*/p^\nu.
$$
\end{pf}

Let $i=n-1$. Since
${\mathrm ord}_pq^n
\geq 2\nu+3$ if $p\geq 3$ and
${\mathrm ord}_pq^n\geq 2\nu+5$ if $p=2$, we
can choose an integer $n_0$ such that $0< n_0 < n-1$ and
$$n_0\cdot{\mathrm ord}_pq>\nu+1/(p-1), \quad
(n-1-n_0)\cdot{\mathrm ord}_pq>\nu+1/(p-1).$$
Due to Lemma \ref{left}, the class
$(\beta_{n-1}\rho_{n-1})(\eta_\E)$ goes to
zero in
$A_{n_0}^*/p^\nu$ via the natural inclusion $A_{n-1}
\to A_{n_0}$. Therefore, by the exact sequence \eqref{ex},
we have a class
$$
\wt{\rho_{n-1}(\eta_\E)}\in
H^1_\et(A_{n_0},\Z/p^\nu(2))
$$
such that $\alpha_{n_0}(
\wt{\rho_{n-1}(\eta_\E)})=\rho_{n-1}(\eta_\E)$.
Next we go on to another inclusion $A_{n_0}
\to A_{0}$.
Write $A_{i,\ol{K}}:=A_{i}\ot_K\ol{K}$.
There is a commutative diagram
\begin{equation}\label{cd}
\begin{CD}
0@.0\\
@VVV@VVV\\
H^1_\et(K,\Z/p^\nu(2))
@>{=}>>
H^1_\et(K,\Z/p^\nu(2))
\\
@VVV@VV{r_1}V\\
H^1_\et(A_{n_0},\Z/p^\nu(2))
@>{\iota_1}>>
H^1_\et(A_{0},\Z/p^\nu(2))
\\
@VVV@VV{r_2}V\\
H^1_\et(A_{n_0,\ol{K}},\Z/p^\nu(2))
@>{\iota_2}>>
H^1_\et(A_{0,\ol{K}},\Z/p^\nu(2)).
\end{CD}
\end{equation}
with exact columns.
\begin{lem}\label{invariancy}
$(r_2\iota_1)(N\cdot\wt{\rho_{n-1}(\eta_\E)})=0$.
In other wards,
$\iota_1(N\cdot\wt{\rho_{n-1}(\eta_\E)})$ is in the image of
$H^1_\et(K,\Z/p^\nu(2))$.
\end{lem}
\begin{pf}
For a finite extension $L/K$, we denote by $R_L$ the integer ring
of $L$ and by $\pi_L$ a
uniformizer of $L$.
Then we have
\begin{align*}
H^1_\et(A_{i,\ol{K}},\Z/p^\nu(2))
&=
H^1_\et(A_{i,\ol{K}},\Z/p^\nu(1))\ot\Z/p^\nu(1)\\
&=
\lim{L/K \text{ finite}}
H^1_\et(A_{i}\ot_KL,\Z/p^\nu(1))\ot\Z/p^\nu(1)\\
&=
\lim{L/K \text{ finite}}
H^1_\et(R_L[[(q^is)]][q^{-1},s^{-1}],\Z/p^\nu(1))\ot\Z/p^\nu(1)\\
&=\lim{L/K \text{ finite}}
R_L[[(q^is)]][q^{-1},s^{-1}]^*/p^\nu\ot\Z/p^\nu(1)
\end{align*}
and
\begin{align*}
R_L[[(q^is)]][q^{-1},s^{-1}]^*&=
R_L[[(q^is)]]^*\times \pi_L^\Z\times s^\Z\\
&=R_L^*\times (1+q^isR_L[[(q^is)]])
\times \pi_L^\Z\times s^\Z.
\end{align*}
Therefore we have
$$
H^1_\et(A_{i,\ol{K}},\Z/p^\nu(2))
=\lim{L/K \text{ finite}}
\left(
(1+q^isR_L[[(q^is)]])\times s^\Z
\right)
/p^\nu\ot\Z/p^\nu(1).
$$
Since $n_0>\nu+1/(p-1)$,
the component $(1+q^{n_0}sR_L[[(q^{n_0}s)]])/p^\nu$ goes to zero
via $\iota_2$:
$$
(1+q^{n_0}sR_L[[(q^{n_0}s)]])/p^\nu
\os{0}{\lra}
(1+sR_L[[s]])/p^\nu.
$$
Therefore, to prove the lemma, it suffices to show that
$\iota_1(\wt{\rho_{n-1}(\eta_\E)})$ goes to a $N$-torsion element
via the composition of the following:
\begin{multline*}
H^1_\et(A_{0},\Z/p^\nu(2))
\os{r_2}{\lra}
H^1_\et(A_{0,\ol{K}},\Z/p^\nu(2))\\
\lra
H^1_\et(\ol{K}((s)),\Z/p^\nu(2))=s^\Z/p^\nu\ot \Z/p^\nu(1).
\end{multline*}
However, the image of $H^1_\et(A_{0},\Z/p^\nu(2))$
must be contained in $G_K$-invariant part of $s^\Z/p^\nu\ot \Z/p^\nu(1)$.
Since it is $N$-torsion, the assertion follows.
\end{pf}
We put by $\hat{\eta}$ the element of $H^1_\et(K,\Z/p^\nu(2))$
such that $r_1(\hat{\eta})=\iota_1(N\cdot\wt{\rho_{n-1}(\eta_\E)})$.
\begin{lem}\label{r1lemma}
$\hat{\eta}
=
N\cdot \wt{\rho(\eta)}$
in $H^1_\et(K,\Z/p^\nu(2))$.
\end{lem}
\begin{pf}
It is enough to show $r_1(\hat{\eta})=
\iota_1(N\cdot\wt{\rho_{n-1}(\eta_\E)})=
r_1(N\cdot \wt{\rho(\eta)})$
in $H^1_\et(A_0,\Z/p^\nu(2))$.
Since they are in the image of
$H^1_\et(K,\Z/p^\nu(2))$,
we may specialize them
via the map $A_{0}\to K$ given by $s\mapsto q$:
$$
\begin{CD}
H^1_\et(K,\Z/p^\nu(2))
\\
@V{r_1}VV\\
H^1_\et(A_0,\Z/p^\nu(2))
@>{x\mapsto x\vert_{s=q}}>>
H^1_\et(K,\Z/p^\nu(2)).
\end{CD}
$$
We want to show that
$\iota_1(N\cdot\wt{\rho_{n-1}(\eta_\E)})\vert_{s=q}=
r_1(N\cdot \wt{\rho(\eta)})\vert_{s=q}=N\cdot \wt{\rho(\eta)}$.
Note $\E_{n-1}\vert_{s=q}=E_{n,K}=K^*/q^{n\Z}$,
$f_\E(v)\vert_{s=q}=f(v)$ and $g_\E(v)\vert_{s=q}=g(v)$ and therefore
$\eta_\E\vert_{s=q}=\eta$.
This shows that the specialization $\iota_1(\rho_{n-1}(\eta_\E))
\vert_{s=q}$ coincides with $\rho(\eta)$
in the cohomology
$H^1_\et(K,H^1_\et(E_{n,\ol{K}},\Z/p^\nu(2)))$.
Hence we have
$\iota_1(\wt{\rho_{n-1}(\eta_\E)})\vert_{s=q}$ coincides with
$\wt{\rho(\eta)}$ modulo $\ker~a$.
Since $\sharp \ker~a\leq N$ the lemma follows.
\end{pf}
Next,
we take another specialization of
$\iota_1(N\cdot\wt{\rho_{n-1}(\eta_\E)})$
via the natural inclusion $\O_0[(qs)^{-1}]\hookrightarrow M:=K((s))$:
$$
\begin{CD}
H^1_\et(K,\Z/p^\nu(2))
\\
@VVV\\
H^1_\et(A_{0},\Z/p^\nu(2))
@>{x\mapsto x\vert_M}>>
H^1_\et(M,\Z/p^\nu(2)).
\end{CD}
$$
Let $E_{n-1,M}=M^*/(q^{n-1}s)^\Z$ be the Tate curve
defined over the field $M$, $\eta_{\E,M}\in
H^0_\Zar(E_{n-1,M},{\cal K}_2)/p^\nu$ the restriction of $\eta_\E$,
and
$$
\rho_{n-1,M}:
H^0_\Zar(E_{n-1,M},{\cal K}_2)/p^\nu \lra
H^1_\et(M,H^1_\et(E_{n-1,\ol{M}},\Z/p^\nu(2)))
$$
the regulator map.
We have the classes
$$
\rho_{n-1,M}(\eta_{\E,M})\in
H^1_\et(M,H^1_\et(E_{n-1,\ol{M}},\Z/p^\nu(2)))
$$
and
$$
\wt{\rho_{n-1,M}(\eta_{\E,M})}\in
H^1_\et(M,\Z/p^\nu(2))
$$
as before.
Then we have
\begin{equation}\label{m1}
r_1(\hat{\eta})\vert_M=\iota_1(N\cdot\wt{\rho_{n-1}(\eta_\E)})\vert_M=
N\cdot\wt{\rho_{n-1,M}(\eta_{\E,M})}\in
H^1_\et(M,\Z/p^\nu(2)).
\end{equation}
We can
think of the right hand side as an element of $H^1_\et(K,\Z/p^\nu(2))$:
\begin{equation}\label{m11}
\hat{\eta}=
N\cdot\wt{\rho_{n-1,M}(\eta_{\E,M})}\in
H^1_\et(K,\Z/p^\nu(2)).
\end{equation}
We calculate the right hand side of \eqref{m11}.
\begin{lem}\label{mlemma}
$
N\cdot\wt{\rho_{n-1,M}(\eta_{\E,M})}=
\pm N\cdot \rho_{C}(\eta_0)
\in H^1_\et(K,\Z/p^\nu(2)).
$
\end{lem}
\begin{pf}
Take a semistable model $X$ of $E_{n-1,M}$ over $K[[s]]$,
and let $Y\subset X$ be the special fiber:
$$
\begin{CD}
Y@>>> X @<<< E_{n-1,M}\\
@VVV@V{\pi_X}VV@VVV\\
\Spec K@>>>\Spec K[[s]] @<<< \Spec M.
\end{CD}
$$
We can easily see that the symbol
$\eta_{\E,M}$ comes from a symbol in the $K$-cohomology
$H^0_\Zar(X,{\cal K}_2)/p^\nu$ of $X$, which we denote by $\eta_X$.
Note that there is a commutative diagram
$$
\begin{CD}
H^0_\Zar(X,{\cal K}_2)/p^\nu@>>>
H^1_\et(K[[s]],R^1\pi_{X*}\Z/p^\nu(2))\\
@VVV@VVV\\
H^0_\Zar(E_{n-1,M},{\cal K}_2)/p^\nu
@>{\rho_{n-1,M}}>>
H^1_\et(M,H^1_\et(E_{n-1,\ol{M}},\Z/p^\nu(2))).
\end{CD}
$$
Moreover, by the proper base change theorem, we have
\begin{align*}
H^1_\et(K[[s]],R^1\pi_{X*}\Z/p^\nu(2))&\os{\cong}{\lra}
H^1_\et(K,(R^1\pi_{X*}\Z/p^\nu(2))\vert_{s=0})\\
&\cong H^1_\et(K,H^1_\et(Y_{\ol{K}},\Z/p^\nu(2)))\\
&\cong H^1_\et(K,\Z/p^\nu(2))\\
&\hookrightarrow H^1_\et(M,\Z/p^\nu(2))
\end{align*}
where the 3rd isomorphism is due to the fact that $Y$ is a chain
of rational curves.
Therefore $\eta_X$ defines a class $\eta'_X\in
H^1_\et(K,\Z/p^\nu(2))$,
and it coincides with the class $\wt{\rho_{n-1,M}(\eta_{\E,M})}$
up to $N$-torsion:
\begin{equation}\label{eta1}
N\cdot \eta'_X=\pm
N\cdot \wt{\rho_{n-1,M}(\eta_{\E,M})}
\in H^1_\et(K,\Z/p^\nu(2)).
\end{equation}
On the other hand, there is also a commutative diagram
$$
\begin{CD}
H^0_\Zar(Y,{\cal K}_2)/p^\nu@>{\rho_{Y}}>>
H^1_\et(K,H^1_\et(Y\ot_K\ol{K},\Z/p^\nu(2)))
@>{\cong}>>
H^1_\et(K,\Z/p^\nu(2))\\
@AAA@AAA\\
H^0_\Zar(X,{\cal K}_2)/p^\nu@>>>
H^1_\et(K[[s]],R^1\pi_{X*}\Z/p^\nu(2)).
\end{CD}
$$
Here $\rho_Y$ is as in \eqref{rhoz}.
Thus we have
\begin{equation}\label{eta2}
\eta'_X=\pm
\rho_Y(\eta_X\vert_Y)
\in H^1_\et(K,\Z/p^\nu(2)).
\end{equation}
By the definition,
\begin{equation}\label{eta3}
\eta_X\vert_Y=
\left\{
\left(\frac{1-\zeta_1^{-1}v}{1-v}\right)^{m_1}
\left(\frac{1-\zeta_1^{-1}\zeta_2}{1-\zeta_2}\right)^{-m_1},
\left(\frac{1-\zeta_2^{-1}v}{1-v}\right)^{m_2}
\left(\frac{1-\zeta_2^{-1}\zeta_1}{1-\zeta_1}\right)^{-m_2}
\right\}
\end{equation}
This means the following (cf. \S \ref{reductionsec}).
Let $Y_0\subset Y$ be the identity
component.
Let $0,~\infty\in Y_0$ be the singular points of $Y$
which are contained in $Y_0$. By attaching $0$ with $\infty$,
we have a nodal curve $C':=Y_0/0\sim \infty$.
There is the morphism ${\mathrm col}:Y \to C'$ which
collapses the other components. Viewing $\eta_0$
as a symbol of $C'$, \eqref{eta3} means
\begin{equation}\label{eta4}
\eta_X\vert_{Y}={\mathrm col}^*\eta_0 \in
H^0_\Zar(Y,{\cal K}_2)/p^\nu.
\end{equation}
\eqref{eta1}, \eqref{eta2} and \eqref{eta4}
yield
$$
N\cdot \wt{\rho_{n-1,M}(\eta_{\E,M})}
=\pm N\cdot \eta'_X=\pm
N\cdot\rho_{Y}(\eta_X\vert_{Y})
=\pm N\cdot\rho_{Y}({\mathrm col}^*\eta_0)
=\pm N\cdot\rho_{C}(\eta_0)
$$
in $H^1_\et(K,\Z/p^\nu(2))$. This completes the proof.
\end{pf}

Now Theorem \ref{partii} is
straightforward from \eqref{m11} and
Lemmas \ref{r1lemma} and \ref{mlemma}.

\subsection{Proof of Part II : Step 2}\label{step4sec}
We finish the proof of (Part II).

Let $\psi_n:E_K\to E_{n,K}$ be the surjective homomorphism given by
$x\mapsto x^n$.
The map $\rho$ in \eqref{rw} and the map
\eqref{qpmap0} are compatible under the pull-back $\psi_n^*$.
Therefore by Lemma \ref{psilem}, we have
\begin{equation}\label{1stk2}
\rho(\eta)\in a(\hsymbol)\subset H^1_\et(K,H^1_\et(E_{\ol{K}},\Z/p^\nu(2))).
\end{equation}
By Theorem \ref{partii}, we have
\begin{equation}\label{2ndk2}
N\cdot a\rho_C(\eta_0)=
\pm N\cdot \rho(\eta)
\in a(\hsymbol).
\end{equation}
Let $\hat{\rho}_C$ be the composition
$$
\hat{\rho}_C:
H^0_\Zar(C,{\cal K}_2) \lra \plim{\nu}H^0_\Zar(C,{\cal K}_2)/p^\nu
\lra
H^1_\et(K,H^1_\et(C_{\ol{K}},\Z_p(2)))
\cong H^1_\et(K,\Z_p(2)).
$$
Then $\rho_C(\eta_0)=\hat{\rho}_C(\eta_0)~{\mathrm mod}~p^\nu$, and
hence we have from \eqref{2ndk2} that
\begin{equation}\label{3rdk2}
\hbox{Image of }N\cdot\hat{\rho}_C(\eta_0)
\in \hsymbol.
\end{equation}
On the other hand, by Corollary \ref{bloch2cor} we have
\begin{equation}\label{4thk2}
\hat{\rho}_C(\eta_0)
=
\pm m_1m_2c_{3,2}([\zeta_1^{-1}\zeta_2]-[\zeta_1^{-1}]-[\zeta_2])
\in H^1_\et(K,\Q_p(2))
\end{equation}
where $c_{3,2}$ is the regulator map
$$
c_{3,2}:B(K)_\Q \cong K_3^{\mathrm ind}(K)_\Q \lra H^1_\et(K,\Q_p(2)).
$$
Put by
$
\hsymbol^{\prime}
$
(resp. $H_{\mathrm ab}^\prime$)
the image of $\hsymbol$ (resp. $H_{\mathrm ab}$)
by the following map
$$
H^1_\et(K,\Z/p^\nu(2))/\ker~a\lra H^1_\et(K,\Q_p/\Z_p(2))/\ker~a.
$$
Since the kernel of $H^1_\et(K,\Z/p^\nu(2))
\to H^1_\et(K,\Q_p/\Z_p(2))$ is dominant by a finite
group $H^0_\et(K,\Q_p/\Z_p(2))$,
to say that $\hsymbol \supset m H_{\mathrm ab}$ for some $m\not=0$
is equivalent to say that
$\hsymbol^{\prime} \supset m' H_{\mathrm ab}^\prime$ for some $m'\not=0$
which does not depend on $\nu$.
Noting the commutative diagram
$$
\begin{CD}
H^1_\et(K,\Z_p(2))
@>{{\mathrm mult.by}~p^{-\nu}}>> H^1_\et(K,\Q_p(2))\\
@VVV@VVV\\
H^1_\et(K,\Z/p^\nu(2))@>>> H^1_\et(K,\Q_p/\Z_p(2))
\end{CD}
$$
we have from \eqref{3rdk2} and \eqref{4thk2} that
\begin{equation}\label{5thk2}
\hbox{Image of }\left(
\frac{Nm_1m_2}{p^\nu}c_{3,2}([\zeta_1^{-1}\zeta_2]-[\zeta_1^{-1}]-[\zeta_2])
\right)\in \hsymbol^\prime.
\end{equation}
Assume $m=m_1=m_2$. Due to \eqref{zetarel} we have
$$
\sum_{\zeta_1}[\zeta_1^{-1}\zeta_2]-[\zeta_1^{-1}]-[\zeta_2]
=-m[\zeta_2]
$$
in the Bloch group $B(K)_\Q$
where $\zeta_1$ runs over all $m$-th roots of unity such that
$\zeta_1\not=\zeta_2$.
Therefore we have from \eqref{5thk2} that
\begin{equation}\label{6thk2}
\hbox{Image of }\left(
\frac{Nm^3}{p^\nu}c_{3,2}[\zeta]
\right)
\in \hsymbol^\prime
\end{equation}
for any $\zeta\in K^*$ such that $\zeta^m=1$.

Let $F\subset K$ be any finite abelian extension over $\Q$.
Then $F$ is contained in a cyclotomic field $\Q(\mu)$.
Since $H^1_\et(\Q(\mu),\Q_p(2))$ is spanned by
$c_{3,2}[\mu^i]$ (cf. Theorem \ref{k3et}
below), so is $H^1_\et(F,\Q_p(2))$.
Therefore, by \eqref{6thk2} and the norm argument
there exists an integer
$m'\not=0$ which does not depend on $\nu$ (but does on $F$) such that
\begin{equation}\label{mot5}
\hbox{Image of }
H^1_\et(F,\Z_p(2))\ot \frac{m'}{p^\nu}\Z/\Z
 \subset \hsymbol^\prime.
\end{equation}
This means $\hsymbol^\prime \supset m'H_{\mathrm ab}^\prime$, and
hence $\hsymbol \supset mH_{\mathrm ab}$ for some $m\not=0$.
This completes the proof of (Part II).

\section{Proof of Theorem \ref{mainthm} : Part III}
\label{galoissec}
In this section, we prove
\begin{description}
\item[(Part III)]
The map
\eqref{kctok}
is surjective if $K\subset \Q_p(\zeta)$ for some root of unity
$\zeta$.
\end{description}
To do this, the following results are crucial.
\begin{thm}\label{k3et}
\begin{enumerate}
\renewcommand{\labelenumi}{(\theenumi)}
\item\label{k3et1} $($\cite{soule1} Thm.1$)$.
Let $F$ be a number field.
Then the regulator map
$$
c_{3,2}:B(F)\ot\Q_p\cong K_3^{\mathrm ind}(F)\ot\Q_p \lra
H_\et^1(F,\Q_p(2))
$$
is bijective.
The dimension of both sides is $r_2$, where $r_2$
denotes the number of
complex places of $F$.
\item\label{k3et12}
Let $\zeta$ be a primitive $n$-th root of unity.
Then the basis of the Bloch group $B(\Q(\zeta))_\Q$ is given by
$\{[\zeta^i]~;~1\leq i < n/2,~(i,n)=1\}$.
$($cf. \cite{irvine} Thm.7.2.4, \cite{neukirch}.$)$
\item\label{k3et2}
$\dim H^1_\et(K,\Q_p(2))=[K:\Q_p]$.
\item\label{k3et3}
If $l\not=p$, then $H^1_\et(K,\Q_l(2))=0$.
\end{enumerate}
\end{thm}
Note that \eqref{k3et2} and \eqref{k3et3} follow from
the Euler-Poincare characteristic (\cite{galois} II 5.7).
Due to \eqref{k3et1} and \eqref{k3et12} one can see that
the map \eqref{kctok} is surjective
if and only if any $x\in H^1_\et(K,\Q_p(2))$ can be written
as a linear combination of
$c_{3,2}([\zeta])$'s.
\subsection{Proof of Part III}
We may assume
$K=\Q_p(\zeta_{m})$
where $m \geq 1$ is an integer and
$\zeta_{m}$ is a primitive $m$-th roots of unity.
In fact, if we show the surjectivity of \eqref{kctok} for
$K=\Q_p(\zeta_m)$,
then we have it for any $K\subset \Q_p(\zeta_m)$
by using the norm map.

We use the following result:
\begin{lem}\label{jan3}
Let $F$ be a number field, and
$P$ the set of all finite places of $F$.
For $v\in P$, we denote by $F_v$ the completion of $F$ by $v$.
Then the natural map
$$
H^1_\et(F,\Q_p(j)) \lra \prod_{v\in P}
H^1_\et(F_v,\Q_p(j))
$$
is injective for $j\not=0$.
\end{lem}
\begin{pf}
See \cite{jangalois} Theorem 3 a).
\end{pf}
Let $l$ be a prime number such that $l\equiv -1$ mod 4,
$(l,m)=1$
and $p$ is complete split
in $\Q(\sqrt{-l})$
(equivalently,
$l\equiv-1$ mod 8 if $p=2$ and $\left(\frac{-l}{p}\right)=1$
if $p\geq3$).
Put $F=\Q(\zeta_m,\sqrt{-l})$.
There are two finite places $v_1$ and $v_2$ of $\Q(\sqrt{-l})$
lying over $p$.
Denote by $e$, $f$ and $g$ the customary meaning for $
F/\Q(\sqrt{-l})$.
Then there are $g$-finite places $\frak{p}_i$
(resp. $\frak{p}'_i$) $1\leq i\leq g$,
of $\Q(\zeta_m,\sqrt{-l})$
lying over $v_1$ (resp. $v_2$).
The completions $F_{\frak{p}_i}$ and
$F_{\frak{p}'_i}$ are isomorphic to $K=\Q_p(\zeta_m)$.
We have $[F:\Q]=2efg=2\varphi(m)$ and $[\Q_p(\zeta_m):\Q_p]=ef$.

Due to Lemma \ref{jan3} and Theorem \ref{k3et} \eqref{k3et3}
the natural map
\begin{equation}\label{fvi1}
H^1_\et(F,\Q_p(2)) \lra \prod_{i=1}^g
H^1_\et(F_{\frak{p}_i},\Q_p(2))\op H^1_\et(F_{\frak{p}'_i},\Q_p(2))
\end{equation}
is injective.
By Theorem \ref{k3et} \eqref{k3et1} and \eqref{k3et2} we have
$\dim H^1_\et(F,\Q_p(2))=efg$ and
$\dim H^1_\et(F_{\frak{p}_i},\Q_p(2))=
\dim H^1_\et(F_{\frak{p}'_i},\Q_p(2))=ef$.
Let
$$
f_1:H^1_\et(F,\Q_p(2)) \lra
\prod_{i=1}^g H^1_\et(F_{\frak{p}_i},\Q_p(2)),\quad
f_2:H^1_\et(F,\Q_p(2)) \lra
\prod_{i=1}^g H^1_\et(F_{\frak{p}'_i},\Q_p(2))
$$
be the natural ones.
We show that $f_1$ and $f_2$
are bijective.
Let $\sigma:F\to F$ be the automorphism such that
$\sigma\sqrt{-l}=-\sqrt{-l}$ and $\sigma(\zeta_m)=\zeta_m$.
It extends to an isomorphism $\bar{\sigma}:\prod_iF_{\frak{p}'_i}\to
\prod_i F_{\frak{p}_i}$
such that the diagram
\begin{equation}\label{fci1}
\begin{CD}
H^1_\et(F,\Q_p(2)) @>{f_1}>>
\prod_iH^1_\et(F_{\frak{p}_i},\Q_p(2))\\
@V{\sigma^*}VV@VV{\bar{\sigma}^*}V\\
H^1_\et(F,\Q_p(2)) @>{f_2}>>
\prod_iH^1_\et(F_{\frak{p}'_i},\Q_p(2))
\end{CD}
\end{equation}
is commutative. On the other hand, let $\tau:F\to F$ be the
automorphism such that $\tau\sqrt{-l}=\sqrt{-l}$ and
$\tau(\zeta_m)=\zeta^{-1}_m$. Since $\tau$ does not
change $v_i$, it extends to an isomorphism
$\bar{\tau}:\prod_i F_{\frak{p}_i}\to
\prod_i F_{\frak{p}_i}$ which makes a commutative
diagram
\begin{equation}\label{fci2}
\begin{CD}
H^1_\et(F,\Q_p(2)) @>{f_1}>>
\prod_iH^1_\et(F_{\frak{p}_i},\Q_p(2))\\
@V{\tau^*}VV@VV{\bar{\tau}^*}V\\
H^1_\et(F,\Q_p(2)) @>{f_1}>>
\prod_iH^1_\et(F_{\frak{p}_i},\Q_p(2)).
\end{CD}
\end{equation}
We see the action of $\sigma^*$ and $\tau^*$ on $H^1_\et(F,\Q_p(2))$
explicitly.
To do this, it is enough to see it on the Bloch group $B(F)_\Q$
by Theorem \ref{k3et} (1).
Let $\zeta_l$ be a primitive $l$-th root of unity.
Since $l\equiv-1$ mod 4 we have $\Q(\sqrt{-l})\subset \Q(\zeta_l)$.
Embedding $B(F)_\Q\hookrightarrow B(\Q(\zeta_m,\zeta_l))_\Q$,
we can see that $B(F)_\Q$ is generated by
$$
\beta_1(\zeta_m^k):=\sum_{i=1}^{l-1}
[\zeta_l^i\zeta_m^k]
,\quad
\beta_2(\zeta_m^k):=\sum_{i=1}^{l-1}
\left(
\frac{i}{l}
\right)
[\zeta_l^i\zeta_m^k]
,\quad (1\leq k <m,~(k,m)=1)
$$
with relations
$$
\beta_1(\zeta_m^k)=-\beta_1(\zeta_m^{-k}),
\quad
\beta_2(\zeta_m^k)=\beta_2(\zeta_m^{-k}).
$$
Here the second equality is due to the fact that
$\left(
\frac{-1}{l}
\right)
=-1$.
Letting $r\in (\Z/l)^*$ be a generator,
we see
\begin{equation}\label{s1ff}
\sigma^*\beta_1(\zeta_m^k)
=
\sum_{i=1}^{l-1}
[\zeta_l^{ir}\zeta_m^k]
=
\beta_1(\zeta_m^k),
\end{equation}
\begin{equation}\label{s2ff}
\sigma^*\beta_2(\zeta_m^k)
=
\sum_{i=1}^{l-1}
\left(
\frac{i}{l}
\right)
[\zeta_l^{ir}\zeta_m^k]
=
-\sum_{i=1}^{l-1}
\left(
\frac{i}{l}
\right)
[\zeta_l^i\zeta_m^k]
=-\beta_2(\zeta_m^k),
\end{equation}
and
\begin{equation}\label{tauff}
\tau^*\beta_1(\zeta_m^k)
=
\beta_1(\zeta_m^{-k})
=-\beta_1(\zeta_m^k),
\quad
\tau^*\beta_2(\zeta_m^k)
=
\beta_2(\zeta_m^{-k})
=\beta_2(\zeta_m^{k}).
\end{equation}

Now we show that $f_i$ are bijective.
Since the dimensions of the target and source are same,
it is enough to show that $f_i$ are injective.
Suppose that
$\sum a_k\beta_1(\zeta_m^k)+b_k\beta_2(\zeta_m^k)$
is in the kernel of $f_1$.
Due to the diagram \eqref{fci2} and \eqref{tauff},
both of $\sum a_k\beta_1(\zeta_m^k)$
and $\sum b_k\beta_2(\zeta_m^k)$ are contained in the kernel of $f_1$.
Then, by the diagram \eqref{fci1},
both of $\sigma^*(\sum a_k\beta_1(\zeta_m^k))
=\sum a_k\beta_1(\zeta_m^k)$
and $\sigma^*(\sum b_k\beta_2(\zeta_m^k))
=-\sum b_k\beta_2(\zeta_m^k)$ are contained in the kernel of $f_2$.
This shows $\ker f_1=\ker f_2$. On the other hand,
since \eqref{fvi1} is injective,
we have $\ker f_1\cap\ker f_2=0$. Thus we have $\ker f_i=0$ for each
$i$.

This completes the proof of (Part III) and hence Theorem \ref{mainthm}.

\begin{rem}
I don't think that the map \eqref{kctok} is always
surjective.
However, I know of no examples where it is not surjective.
\end{rem}

\section{Applications of Theorem \ref{mainthm}}\label{applysec}
For an abelian group $M$ we denote by $M_{\mathrm tor}$ and
$M_{\mathrm div}$ the torsion subgroup and the maximal divisible
subgroup of $M$
respectively:
$$
M_{\mathrm tor}:=\bigcup_{m\geq1}M[m],\quad
M_{\mathrm div}:=\bigcap_{m\geq1}mM.
$$
\subsection{Consequence of Suslin's exact sequence}
Let $F$ be any field, and $X$ a nonsingular curve over $F$.
Let $p$ be a prime number such that
$p\not={\mathrm char}(F)$.
Passing to the inductive limit over $p^\nu$, we have
from \eqref{universalc}
\begin{equation}\label{sus}
0\lra H^0_\Zar(X,{\cal K}_2)\ot\Q_p/\Z_p
\lra H^2_\et(X,\Q_p/\Z_p(2)) \lra
H^1_\Zar(X,{\cal K}_2)[p^\infty]\lra 0.
\end{equation}
The following is an easy consequence of Suslin's exact sequence.
\begin{prop}\label{equivthm}
Assume that for any finite $p$-torsion $G_F$-module $M$,
$H^i_\et(F, M)$ is finite for all $i\geq 0$.
Then the following are equivalent.
\begin{enumerate}
\renewcommand{\labelenumi}{(\theenumi)}
\item\label{0}
The map $H^0_\Zar(X,{\cal K}_2)\ot\Z_p
\to H^2_\et(X,\Z_p(2))$ is surjective.
\item\label{1}
The map $H^0_\Zar(X,{\cal K}_2)\ot\Q_p
\to H^2_\et(X,\Q_p(2))$ is surjective.
\item\label{2}
The map $H^0_\Zar(X,{\cal K}_2)\ot\Q_p/\Z_p
\to H^2_\et(X,\Z_p(2))\ot\Q_p/\Z_p$ is surjective
$($and hence bijective by \eqref{sus}$)$.
\item\label{3}
The corank of $H^0_\Zar(X,{\cal K}_2)\ot\Q_p/\Z_p$ is greater than or
equal to the corank of
$H^2_\et(X,\Q_p/\Z_p(2))$.
\item\label{4}
$H^1_\Zar(X,{\cal K}_2)[p^\infty]\cong H_\et^{3}(X,\Z_p(2))[p^\infty]$.
\item\label{5}
$H^1_\Zar(X,{\cal K}_2)[p^\infty]$ is finite.
\item\label{66}
$H^0_\Zar(X,{\cal K}_2)/p^\nu
\cong H^2_\et(X,\Z_p(2))/p^\nu$ for all $\nu\geq 1$.
\item\label{7}
$H^1_\Zar(X,{\cal K}_2)[p^\nu]\cong H_\et^{3}(X,\Z_p(2))[p^\nu]$
for all $\nu\geq 1$.
\end{enumerate}
\end{prop}
\begin{pf}
By the Hochschild-Serre spectral sequence,
the assumption on $F$ implies that $H^i_\et(X,\Z/p^\nu(j))$ is finite
for all $i,~j$ and $\nu$,
and hence we have that
$H^i_\et(X,\Z_p(j))$ is a finitely generated $\Z_p$-module
and there are exact sequences
\begin{equation}\label{milne165}
0\lra H_\et^i(X,\Z_p(j))/p^\nu \lra H_\et^i(X,\Z/p^\nu(j))
\lra H_\et^{i+1}(X,\Z_p(j))[p^\nu]\lra0
\end{equation}
for all $i,~j$ and $\nu$ (see for example \cite{milne} p.165, Lemma 1.11).

\medskip

Since $H^2_\et(X,\Z_p(2))$ is a finitely generated $\Z_p$-module,
\eqref{1} is equivalent to \eqref{2}.
Due to the exact sequences \eqref{sus}
and \eqref{milne165},
we have a commutative diagram
$$
\begin{CD}
0@>>> H^0_\Zar(X,{\cal K}_2)\ot\Q_p/\Z_p
@>>> H^2_\et(X,\Q_p/\Z_p(2))
@>>>
H^1_\Zar(X,{\cal K}_2)[p^\infty]@>>> 0\\
@.@VVV@VV{=}V@VVV\\
0@>>> H^2_\et(X,\Z_p(2))\ot\Q_p/\Z_p
@>>> H^2_\et(X,\Q_p/\Z_p(2))
@>>> H^3_\et(X,\Z_p(2))[p^\infty]@>>> 0.
\end{CD}
$$
Therefore, \eqref{2} is equivalent to \eqref{4}.
Since the map $H^0_\Zar(X,{\cal K}_2)\ot\Q_p/\Z_p
\ra H^2_\et(X,\Z_p(2))\ot\Q_p/\Z_p$ is always injective,
\eqref{2} is equivalent to \eqref{3}.
Thus we have completed
\begin{center}
\eqref{1}
$\Longleftrightarrow$\eqref{2}$\Longleftrightarrow$\eqref{3}
$\Longleftrightarrow$\eqref{4}.
\end{center}

\smallskip

\noindent\eqref{4}$\Longrightarrow$\eqref{5}.
This follows from the fact that $H^3_\et(X,\Z_p(2))$ is a finitely
generated $\Z_p$-module.

\noindent\eqref{5}$\Longrightarrow$\eqref{3}.
This follows from \eqref{sus}.

\noindent\eqref{2}$\Longrightarrow$\eqref{66} and \eqref{7}.
Due to Suslin's exact sequences \eqref{universalc} and \eqref{sus},
we have a commutative diagram
$$
\begin{CD}
0@>>> H^2_\et(X,\Z_p(2))/p^\nu
@>>> H^2_\et(X,\Z/p^\nu(2))
@>>>
H^3_\et(X,\Z_p(2))[p^\nu]@>>> 0\\
@.@A{a_1}AA@AA{=}A@AA{a_2}A\\
0@>>> H^0_\Zar(X,{\cal K}_2)/p^\nu
@>>> H^2_\et(X,\Z/p^\nu(2))
@>{\partial}>>
H^1_\Zar(X,{\cal K}_2)[p^\nu]@>>> 0\\
@.@VVV@VVV@VV{b}V\\
0@>>> H^0_\Zar(X,{\cal K}_2)\ot\Q_p/\Z_p
@>>> H^2_\et(X,\Q_p/\Z_p(2))
@>>>
H^1_\Zar(X,{\cal K}_2)[p^\infty]@>>> 0.
\end{CD}
$$
We show that $a_1$ is surjective. It implies
that $a_1$ and $a_2$ are bijective.
To do this, it is enough to see that $H^2_\et(X,\Z_p(2))/p^\nu $
goes to zero by the map $\partial$. Since $b$ is injective,
it is enough to see that it goes to zero by $b\partial$.
However since $H^0_\Zar(X,{\cal K}_2)\ot\Q_p/\Z_p =
H^2_\et(X,\Z_p(2))\ot\Q_p/\Z_p$, it is clear.

\noindent\eqref{66}$\Longrightarrow$\eqref{0}.
Nakayama's lemma.

\noindent\eqref{0}$\Longrightarrow$\eqref{1}.
Clear.

\noindent\eqref{66}$\Longrightarrow$\eqref{2}.
Clear.

\noindent\eqref{7}$\Longrightarrow$\eqref{4}.
Clear.
\end{pf}
\begin{prop}\label{equivthm2}
If the equivalent conditions in Proposition \ref{equivthm}
are satisfied, then we have 
$$
\plim{\nu}H^1_\Zar(X,{\cal K}_2)[p^\nu]=0.
$$
In particular, any $p$-divisible subgroup $D\subset
H^1_\Zar(X,{\cal K}_2)$ (i.e. $D=pD$) is uniquely $p$-divisible.
\end{prop}
\begin{pf}
Since \eqref{universalc} is an exact sequence of finite groups,
the exactness is preserved after taking the projective limit:
$$
0\lra \plim{\nu}H^0_\Zar(X,{\cal K}_2)/p^\nu
\lra H^2_\et(X,\Z_p(2))\lra 
\plim{\nu}H^1_\Zar(X,{\cal K}_2)[p^\nu]\lra 0.
$$
The vanishing of the last term follows from (1).
If $D\subset H^1_\Zar(X,{\cal K}_2)$ is a $p$-divisible subgroup,
then we have $\plim{\nu}D[p^\nu]=0$ and
therefore $D$ has no $p$-torsion. 
\end{pf}


\subsection{Applications to $V(E_K)$}\label{torsionve}
Let us go back to the Tate curve $E_K=K^*/q^\Z$.
Suppose that $K\subset \Q_p(\zeta)$ for some root of unity $\zeta$.
Then by Theorem \ref{mainthm}
and Proposition \ref{equivthm}
we have
$$H^1_\Zar(E_K,{\cal K}_2)[p^\nu]\cong H_\et^{3}(E_K,\Z_p(2))[p^\nu],
\quad
V(E_K)[p^\nu]\cong H^2_\et(K,H_\et^1(E_{\ol{K}},\Z_p(2)))[p^\nu]$$
for all $\nu\geq 1$.
Here we note
$$
H^i_\et(K,H_\et^j(E_{\ol{K}},\Z_p(r)))
\os{\mathrm def}{=}\plim{\nu}
H^i_\et(K,H_\et^j(E_{\ol{K}},\Z/p^\nu(r))).
$$
Due to T.Sato \cite{TSato} (cf. Corollary \ref{bysator}),
the above is also true for the $l$-torsion parts.
Thus we have
\begin{equation}\label{vektor}
V(E_K)_{\mathrm tor}=\bigoplus_lV(E_K)[l^\infty]\cong \bigoplus_l
H^2_\et(K,H_\et^1(E_{\ol{K}},\Z_l(2)))[l^\infty].
\end{equation}
\begin{lem}\label{llll}
Let $\nu_0\geq 0$ be the largest integer such that a primitive
$l^{\nu_0}$-th root of unity is contained in $K^*$.
Then there is the natural isomorphism
$$H^2_\et(K,H_\et^1(E_{\ol{K}},\Z_l(2)))[l^\infty]
\cong K_2^M(K)/(l^\nu K_2^M(K)+\{q,K^*\})
$$
for $\nu\geq \nu_0$.
\end{lem}
\begin{pf}
By the weight exact sequence \eqref{van0}, we have an exact sequence
$$
 H^1_\et(K,\Z_l(1))\os{\delta}{\lra} H^2_\et(K,\Z_l(2))
\lra
H^2_\et(K,H^1_\et(E_{\ol{K}},\Z_l(2)))
\lra H^2_\et(K,\Z_l(1)).
$$
Recall the isomorphism
$$
K^M_2(K)/l^{\nu}\os{\cong}{\lra} \plim{i} K^M_2(K)/l^i
\os{\cong}{\lra}H^2_\et(K,\Z_l(2))
$$
for $\nu\geq \nu_0$.
Under this isomorphism and $H^1_\et(K,\Z/l^i(1))\cong K^*/l^i$,
the map $\delta$ is given by $x\mapsto \{x,q\}$.
Therefore the cokernel of
$\delta$ is isomorphic to $
K_2^M(K)/(l^{\nu} K_2^M(K)+\{q,K^*\})$.
On the other hand, $H^2_\et(K,\Z_l(1))$ is isomorphic to $\Z_l$
which has no torsion.
Thus we have the assertion.
\end{pf}
By Lemma \ref{llll} and the isomorphism \eqref{vektor}, we have
\begin{cor}\label{torstr}
Let the notations and assumption be as in Theorem \ref{mainthm}.
Write by $\mu_n$ the group of all roots of unity in $K$ where
$n$ denotes the cardinality.
Then there is the natural isomorphism
\begin{equation}\label{torisom}
V(E_K)_{\mathrm tor}\cong K^M_2(K)/(nK_2^M(K)+\{q,K^*\})
\os{\cong}{\lra}
\mu_n/(q,K^*)_n
\end{equation}
where the last isomorphism is the map
induced from the Hilbert symbol
$(~,~)_{n}:
K^*/n\times K^*/n\to \mu_n$
(cf. \cite{fesenko} IX (4.3)).
In particular, $V(E_K)_{\mathrm tor}$ and hence
$K_1(E_K)_{\mathrm tor}$ are finite.
\end{cor}

There is the exact sequence
$$
0\lra H^1_\Zar(E_K,{\cal K}_2)/m\lra H^3_\et(E_K,\Z/m(2))
\lra H^3_\et(K(E_K),\Z/m(2))
$$
for all $m\geq1$ (\cite{ms} \S 18). From this we have an
isomorphism
\begin{equation}\label{modmisom}
V(E_K)/m\cong K^M_2(K)/(mK_2^M(K)+\{q,K^*\}).
\end{equation}
\begin{cor}\label{torstrdec}
Let the notations and assumption be as in Theorem \ref{mainthm}.
Then we have a decomposition
\begin{equation}\label{decisom}
V(E_K)= V(E_K)_{\mathrm tor}\op V(E_K)_{\mathrm div},
\end{equation}
and $V(E_K)_{\mathrm div}$ is
uniquely divisible.
\end{cor}
\begin{pf}
Due to
\eqref{torisom} and \eqref{modmisom} we have $V(E_K)_{\mathrm tor}/m
\os{\cong}{\to}V(E_K)/m$ for all $m\geq1$. Therefore
$V(E_K)/V(E_K)_{\mathrm tor}$ is uniquely divisible.
Since $V(E_K)_{\mathrm tor}$ is finite whose order is divided by $n$,
we have that
$nV(E_K)=V(E_K)_{\mathrm div}$ and hence it maps onto
$V(E_K)/V(E_K)_{\mathrm tor}$. Moreover $
V(E_K)_{\mathrm div}\cap V(E_K)_{\mathrm tor}=0$ because of the finiteness
of $V(E_K)_{\mathrm tor}$. Thus we have
$V(E_K)_{\mathrm div}\os{\cong}{\to}V(E_K)/V(E_K)_{\mathrm tor}$.
\end{pf}
\begin{rem}
It is known that $\dim_\Q V(E_K)_{\mathrm div}=+\infty$
(\cite{szam} Appendix).
\end{rem}

\subsection{Earlier works on $V(X)$}\label{others}
Several people studied $V(X)$ (mainly from the viewpoint
of the class field theory) and obtained related results
to \S \ref{torsionve}.
Here are some:
\begin{enumerate}
\renewcommand{\labelenumi}{(\theenumi)}
\item (T.Sato \cite{TSato}).
If $X$ is an elliptic curve over a $p$-adic field with bad reduction,
$V(X)_{\mathrm div}$ is uniquely $l$-divisible for $l\not=p$.
\item (S.Saito \cite{ssaito}).
Let $X$ be a nonsingular projective curve over a $p$-adic field.
Then $V(X)/V(X)_{\mathrm div}$ is finite.
\item (Colliot-Th\'el\`ene and Raskind \cite{col-ras}).
If $X$ is a nonsingular projective curve over a $p$-adic field
with good reduction, 
$V(X)$ is a direct sum of
a uniquely divisible group and $V(X)_{\mathrm tor}$.
Moreover, $V(X)[l^\infty]\cong J(k)[l^\infty]$ for any $l\not=p$
where $k$ is the residue field of $K$ and
$J/k$ is the Jacobian variety of the special fiber.
\end{enumerate}
(1) is a consequence of the surjectivity
of the $l$-adic regulator on
$K_2(E_K)
$ (Proposition \ref{equivthm2}).
Suppose that $X$ has a good reduction.
Then
it follows from the Euler-Poincare characteristic (\cite{galois} II 5.7)
that we have $H^2_\et(X,\Q_l(2))=0$.
Therefore $V(X)[l^\infty]$ is finite and 
$V(X)_{\mathrm div}$ is uniquely $l$-divisible for $l\not=p$
(Propositions \ref{equivthm}, \ref{equivthm2}).
However I do not know any previous results about
finiteness of $V(X)[p^\infty]$.
(Note that (2) and (3) do not imply anything about finiteness of
$p$-power torsion.)
When $X$ has a bad reduction and
the genus of $X$ is greater than 1, the question of the finiteness
remains open even
for the $l$-power torsion part.

\subsection{Other Corollaries}
\begin{cor}\label{pnuk2}
Let the notations and assumption be as in Theorem \ref{mainthm}.
Then the $p$-adic regulator $H^0_\Zar(E_K,{\cal K}_2)\ot\Z_p
\to H^2_\et(E_K,\Z_p(2))$ is surjective, and it induces an isomorphism
$$
H^0_\Zar(E_K,{\cal K}_2)/p^\nu
\os{\cong}{\lra} H^2_\et(E_K,\Z_p(2))/p^\nu$$
for all $\nu\geq 1$.
\end{cor}
\begin{pf}
Straightforward from Theorem \ref{mainthm} and
Proposition \ref{equivthm}.
\end{pf}

\begin{cor}
Let the notations and assumption be as in Theorem \ref{mainthm}.
Denote by $n_0$ the cardinality of the subgroup $(q,K^*)_n\subset \mu_n$.
Let $m\geq1$ be an integer.
Put by $e_m$ the order of the kernel of a map
\begin{equation}\label{exactscon}
\Z/(n_0,m) \lra \Z/(n,m),\quad k\longmapsto k\cdot nn^{-1}_0.
\end{equation}
Then the cohomology of the sequence
\begin{equation}\label{exacts}
H^0_\Zar(E_K,{\cal K}_2)/m \os{\tau^\et_\infty}{\lra} K^*/m
\os{x\mapsto\{x,q\}}{\lra}
K_2(K)/m
\end{equation}
at the middle term
is a cyclic group of order $e_m$.
In particular, \eqref{exacts}
is exact if and only if \eqref{exactscon} is injective.
\end{cor}
\begin{pf}
Since we have the exact sequence
$$
\begin{matrix}
&&H^0_\Zar(E_K,{\cal K}_2)/m
\\
&&\downarrow\\
H^1_\et(K,\Z/m(2))&\os{a}{\lra}&
H^1_\et(K,H^1(E_{\ol{K}},\Z/m(2)))
&\lra&
K^*/m
&\os{x\mapsto\{x,q\}}{\lra} &
K_2(K)/m,
\end{matrix}
$$
the cohomology of \eqref{exacts} is isomorphic to
the cokernel of the map
$$
H^0_\Zar(E_K,{\cal K}_2)/m\lra \Coker ~a.
$$
By Proposition \ref{equivthm} and Theorem \ref{mainthm},
we have an exact sequence
$$
0\ra H^0_\Zar(E_K,{\cal K}_2)/m
\ra
H^1_\et(K,H^1(E_{\ol{K}},\Z/m(2)))
\os{b}{\ra}
H^2_\et(K,H^1(E_{\ol{K}},\hat{\Z}(2)))[m]
\ra 0.
$$
Moreover there is a commutative diagram
$$
\begin{CD}
H^1_\et(K,H^1(E_{\ol{K}},\Z/m(2)))
@>{b}>>
H^2_\et(K,H^1(E_{\ol{K}},\hat{\Z}(2)))[m]\\
@A{a}AA@AA{a'}A\\
H^1_\et(K,\Z/m(2))
@>{b'}>>
H^2_\et(K,\hat{\Z}(2))[m]
\end{CD}
$$
where $b'$ is surjective.
Therefore the cohomology of \eqref{exacts} is isomorphic to
the cokernel of $a'$.
As we have seen in
the proof of Lemma \ref{llll}, we have
$$
H^2_\et(K,H^1(E_{\ol{K}},\hat{\Z}(2)))_{\mathrm tor}
\cong \mu_n/(q,K^*)_n\cong\Z/nn_0^{-1},\quad
H^2_\et(K,\hat{\Z}(2))_{\mathrm tor}
\cong \mu_n\cong\Z/n,
$$
and the map $a'$ can be identified with the natural map
$$
\Z/n[m] \lra \Z/nn^{-1}_0[m].
$$
Its cokernel is isomorphic to the kernel of \eqref{exactscon}.
\end{pf}


\section{Surjectivity of $l$-adic regulator on $K_2$ of open Tate curves}
\label{opentatesec}

\begin{thm}\label{gensa}
Let $l\not=p$ be a prime number.
Let $U_K\subset E_K$ be an arbitrary Zariski open set
$($no assumption on $K$$)$.
Then
the $l$-adic regulator
\begin{equation}\label{reg101}
K_2(U_K)\ot\Q_l\lra H^2_\et(U_K,\Q_l(2))
\end{equation} is surjective.
\end{thm}
This is a generalization of the main result of T.Sato's
thesis \cite{TSato} which proved the above in case
$U_K=E_K$.
When $l=p$, the question of the surjectivity for $U_K\not=E_K$
remains open.
\begin{pf}
Using the norm map, we can replace $K$ with an arbitrary finite
extension $L$ over $K$. Thus we may assume that $E_K$ is defined
by an equation
\begin{equation}\label{wei}
y^2=x^3+x^2+c
\end{equation}
with ${\mathrm ord}_K(c)=n\geq1$, and
$D_K:=E_K-U_K=P_1+\cdots+P_s$ with each $P_i\in E_K(K)$.
In the same way as the proof of \eqref{open1}, we have
\begin{equation}\label{open2}
H_\et^{2}(U_K,\Q_l(2))\cong H^1_\et(K,H^1(U_{\ol{K}},\Q_l(2))).
\end{equation}
It follows from the exact sequence (cf. \S \ref{partipre})
\begin{equation}\label{e1}
0\lra \Q_l(2)
\lra H^1_\et(U_{\ol{K}},\Q_l(2))
\lra
\Q_l(1)\op\bigoplus_{i=2}^s\Q_l(1)([P_i]-[P_1])\lra 0
\end{equation}
and Theorem \ref{k3et} \eqref{k3et3} that we have
\begin{align*}
H^1_\et(K,H^1(U_{\ol{K}},\Q_l(2)))&\cong
H^1_\et(K,\Q_l(1))\op\bigoplus_{i=2}^sH^1_\et(K,\Q_l(1))([P_i]-[P_1])\\
&\cong
\Q_l\op\bigoplus_{i=2}^s\Q_l([P_i]-[P_1]).
\end{align*}
Thus we can rewrite the $l$-adic regulator \eqref{reg101}
in the following form:
\begin{equation}\label{open3}
K_2(U_K)\ot\Q_l\lra \Q_l\op\bigoplus_{i=2}^s\Q_l([P_i]-[P_1]).
\end{equation}
Since $K_2(E_K)\ot\Q_l$ is onto the first component
(Corollary \ref{bysator}), it is enough to show that
the composition
\begin{equation}\label{open4}
K_2(U_K)\ot\Q_l\lra \Q_l\op\bigoplus_{i=2}^s\Q_l([P_i]-[P_1])
\os{\mathrm pr}{\lra}\bigoplus_{i=2}^s\Q_l([P_i]-[P_1])
\end{equation}
is surjective.
To do this, we note that the map \eqref{open4} is obtained from
the tame symbols.
More precisely let $\tau_{P_i}$ be the tame symbol at $P_i$
and $\ord_K:K^*\to \Z$ the valuation such that $\ord_K(\pi_K)=1$.
Then \eqref{open4} is obtained by tensoring
\begin{equation}\label{tame31}
\sum_{i=2}^s{\mathrm ord}_K\cdot\tau_{P_i}:
H^0_\Zar(U_K,{\cal K}_2)
\lra
\bigoplus_{i=2}^s\Z([P_i]-[P_1])
\end{equation}
with $\Q_l$.
Therefore it suffices to show that the cokernel of \eqref{tame31}
is finite.

Obviously we can reduce it to the case $s=2$. Moreover, by using
the translation $x\mapsto x-a$, we may assume $P_1=O$. Put
$P=P_2$. Write $P$ as $P=\alpha q^{r/n}$ ($\alpha\in
R^*$) under the identification $K^*/q^\Z=E_K(K)$. Put $Q=\zeta
q^{r/n}$ with $\zeta^m=1$. Then the cokernel of the tame symbol
$\tau_Q:H^0_\Zar(E_K-\{Q,O\},{\cal K}_2)\to K^*$ is finite. In
fact, putting
$f(u)=\theta(\zeta^{-1}q^{-r/n}u)^{nm}/\theta(u)^{nm-1}\theta(q^{-rm}u)$,
the symbol $\{a,f(u)\}$ goes to $a^{nm}$. Therefore, in order to
show that the cokernel of \eqref{tame31} is finite in case
$D_K=P+O$, it suffices to show it in case $D_K=P+Q+O$. It is also
reduced to the case $D_K=P+Q$. By the translation, it is reduced
to the case $D_K=P'+O$ where $P'=\alpha\zeta^{-1}\in R^*$. By
choosing a suitable $\zeta$, we can assume
$\alpha\zeta^{-1}~{\mathrm mod}\pi_K\not\equiv 1$. Moreover,
replacing $K$ with $K(\sqrt{\alpha\zeta^{-1}})$, we may assume
that there is $\beta\in R^*$ such that $\alpha\zeta^{-1}=\beta^2$.
By using the translation, we can reduce the case
$D_K=P^{\prime\prime}+(-P^{\prime\prime})$ with
$P^{\prime\prime}=\beta$ and $(-P^{\prime\prime})=\beta^{-1}$.
Summarizing the above, we have reduced the proof to the
following claim:

\begin{claim}
Suppose $P=\alpha\in R^*$ and
$\alpha~{\mathrm mod}\pi_K\not\equiv \pm 1$.
Then the cokernel of
${\mathrm ord}_K\cdot\tau_P:
H^0_\Zar(E_K-\{P,(-P)\},{\cal K}_2)\to \Z$ is finite.
\end{claim}

Let us go back to the equation \eqref{wei}.
Let $P=(a,b)$ be the coordinate expression by $(x,y)$ with $a,~b\in R$.
Note $(-P)=(a,-b)$.
We consider a $K_2$-symbol
\begin{equation}\label{my2}
\xi:=\left\{\frac{y-\sqrt{a+1}x}{y+\sqrt{a+1}x},
\frac{-c}{x^2(x-a)}\right\}.
\end{equation}
We have
$$
\tau_Q(\xi)=
\begin{cases}
(b-\sqrt{a+1}a)/(b+\sqrt{a+1}a)& Q=P\\
(b+\sqrt{a+1}a)/(b-\sqrt{a+1}a)& Q=(-P)\\
1& \text{otherwise}.
\end{cases}
$$
This shows $\xi\in H^0_\Zar(E_K-\{P,(-P)\},{\cal K}_2)$. The
assumption $\alpha\in R^*$ and $\alpha~{\mathrm
mod}\pi_K\not\equiv \pm 1$ implies $a~{\mathrm mod}\pi_K\not\equiv
0,~-1$. Suppose $p\not=2$. Since
$(b+\sqrt{a+1}a)-(b-\sqrt{a+1}a)=2\sqrt{a+1}a$ is a unit, either
$(b-\sqrt{a+1}a)$ or $(b+\sqrt{a+1}a)$ is a unit in $R$. On the
other hand, the order of $(b-\sqrt{a+1}a)(b+\sqrt{a+1}a)=c$ is
$n$. We have
$$
{\mathrm ord}_K\cdot\tau_P(\xi)={\mathrm ord}_K\left(
\frac{b-\sqrt{a+1}a}{b+\sqrt{a+1}a}\right)=\pm n.
$$
This means that the cokernel of ${\mathrm ord}_K\cdot\tau_P$
is finite in case $p\not=2$.
Suppose $p=2$. If $n/2>{\mathrm ord}_K(2)$,
one can show that the order of $(b-\sqrt{a+1}a)/(b+\sqrt{a+1}a)$
is not zero by the same argument as above.
For a small $n$, we take a finite covering
$E_{m,K}=K^*/q^{m\Z}\to E_K$ given by $x\mapsto x^m$ with $m\gg 1$.
Since we have obtained the surjectivity of the $l$-adic regulator
for any Zariski open set of $E_{m,K}$,
we can obtain it for $E_K$ by using the transfer map for
$E_{m,K}\to E_K$.
This completes the proof for $p=2$ and all $n\geq 1$.
\end{pf}
\begin{cor}
The $l$-power torsion $K_1(U_K)[l^\infty]$ is finite.
\end{cor}
\begin{pf}
This follows from Theorem \ref{gensa} and Proposition \ref{equivthm}.
\end{pf}

\medskip

\begin{flushleft}

\noindent

Graduate School of Mathematics, Kyushu University,
Hakozaki Higashi-ku Hukuoka 812-8581,
JAPAN

\medskip

\noindent
E-mail address : asakura@@math.kyushu-u.ac.jp
\end{flushleft}

\end{document}